\documentclass[11pt]{amsart}
\usepackage{pdfsync}
\usepackage{a4,fullpage}
\usepackage{graphicx}
\usepackage{color}
\usepackage{amsfonts}
\usepackage{enumerate}
\usepackage{latexsym}
\usepackage{epsfig}
\usepackage{amsthm}
\usepackage{amsmath}
\usepackage{amssymb}
\usepackage{latexsym}
\usepackage{epsfig}
\usepackage{amssymb,latexsym}
\usepackage[all]{xy}
\usepackage[ps2pdf, colorlinks=true]{hyperref}
\usepackage{pdfsync}
\usepackage{a4,fullpage}
\usepackage{enumerate}
\usepackage{graphicx}
\usepackage{latexsym}
\usepackage{amsfonts}
\usepackage[utf8]{inputenc}
\usepackage{amsthm}
\usepackage{amsmath}
\usepackage{amssymb}
\newtheorem{theorem}{Theorem}[section]

\newtheorem{conj}[theorem]{Conjecture}

\newtheorem{problem}[theorem]{Problem}

\newcommand{\ii}{\op{i}}

\newcommand{\R}{{\mathbb R}}

\newcommand{\Z}{{\mathbb Z}}

\newcommand{\T}{{\mathbb T}}

\newcommand{\phy}{\varphi}

\newcommand{\op}[1]{\!\!\mathop{\rm ~#1}\nolimits}
\newcommand{\DD}{\!\mathop{\rm d\!}\nolimits}

\newcommand{\scriptop}[1]{\!\!\mathop{\mbox{\rm \scriptsize ~#1}}\nolimits}

\newenvironment{remark}{\refstepcounter{theorem}\par\medskip\noindent{\bf
Remark~\thetheorem~~}}{\unskip\nobreak\hfill\hbox{ $\oslash$}\par\bigskip}

\newenvironment{example}{\refstepcounter{theorem}\par\medskip\noindent{\bf
Example~\thetheorem~~}}{\unskip\nobreak\hfill\hbox{ $\oslash$}\par\bigskip}

\newenvironment{definition}{\refstepcounter{theorem}\par\medskip\noindent{\bf
Definition~\thetheorem~~}}
\newcommand{\got}[1]{\mathfrak{#1}}

\newcommand{\formel}[1]{[\![#1]\!]}

\newcommand{\re}{\mathop{\mathfrak{R}}\nolimits}

\newcommand{\RM}{\mathbb{R}}
\newcommand{\ZM}{\mathbb{Z}}

\newcommand{\h}{\hbar}

\newcommand{\aangles}{Liou\-ville-Arnold-Mi\-neur}

\title{Symplectic theory \\of completely integrable Hamiltonian systems}

\date{}

\author{{\'A}lvaro Pelayo and San V\~u Ng\d oc \\
{\tiny In memory of  Professor Johannes (Hans) J. Duistermaat (1942-2010)}}

\begin{document}

\maketitle

\begin{abstract}
  This paper explains the recent developments on the symplectic theory
  of Hamiltonian completely integrable systems on symplectic
  $4$\--manifolds, compact or not.  One fundamental ingredient of
  these developments has been the understanding of singular affine
  structures. These developments make use of results obtained by many
  authors in the second half of the twentieth century, notably Arnold,
  Duistermaat and Eliasson, of which we also give a concise survey.
  As a motivation, we present a collection of remarkable results
  proven in the early and mid 1980s in the theory of Hamiltonian Lie
  group actions by Atiyah, Guillemin\--Sternberg and Delzant among
  others, and which inspired many people, including the authors, to
  work on more general Hamiltonian systems. The paper
  concludes discussing a spectral conjecture for quantum integrable
  systems.
\end{abstract}

\section{Introduction} \label{sec:intro}

In the mathematical theory of conservative systems of differential
equations one finds cases that are solvable in some sense, or
\emph{integrable}, which enables one to study their dynamical behavior
using differential geometric and Lie\footnote{Sophus Lie has been one
  of the most influential figures in differential geometry.  Many
  modern notions of differential geometry were known to Lie in some
  form, including the notion of symplectic manifold, symplectic and
  Hamiltonian vector fields, transformation (=Lie) groups, and (in a
  particular instance) momentum maps.}  theoretic methods, in
particular the theory of Lie group actions on symplectic manifolds.
Integrable systems are a fundamental class of ``explicitly solvable''
dynamical systems of current interest in differential and algebraic
geometry, representation theory, analysis and physics.  Their study
usually combines ideas from several areas of mathematics, notably
partial differential equations, microlocal analysis, Lie theory,
symplectic geometry and representation theory.  In this paper we focus
on finite dimensional completely integrable Hamiltonian systems
(sometimes called ``Liouville integrable systems'') in the context of
symplectic geometry.

Many authors have studied dynamical problems for centuries.  Galileo
made great advances to the subject in the late XVI and early XVII
centuries, and formulated the ``laws of falling bodies''.  An
important contribution was made by Huygens in the XVII century, who
studied in detail the spherical pendulum, a simple but fundamental
example. Galileo's ideas were generalized and reformulated by William
Hamilton (1805-1865) using symplectic geometry, who said: ``\emph{the
  theoretical development of the laws of motion of bodies is a problem
  of such interest and importance that it has engaged the attention of
  all the eminent mathematicians since the invention of the dynamics
  as a mathematical science by Galileo, and especially since the
  wonderful extension which was given to that science by Newton}''
(1834, cf. J.R. Taylor \cite[p. 237]{taylor}). Many of the modern
notions in the mathematical theory of dynamical systems date back to
the late XIX century and the XX century, to the works of Poincar{\'e},
Lyapunov, Birkhoff, Siegel and the Russian school in the qualitative
theory of ordinary differential equations.

A \emph{completely integrable Hamiltonian system} may be given by the
following data\footnote{We will explain this definition in detail
  later, starting with the most basic notions.}: (1) a
$2n$\--dimensional smooth manifold $M$ equipped with a symplectic
form, and (2) $n$ smooth functions $f_1,\ldots,f_n \colon M \to
\mathbb{R}$ which generate vector fields that are pairwise linearly
independent at almost every point, and which Poisson commute. In local
symplectic coordinates, this commuting condition amounts to the
vanishing of partial differential equations involving the $f_i$,
eg. see Section \ref{sec:semitoric}.  Many times we will omit the word
``Hamiltonian'', and refer simply to ``completely integrable
systems''.  A completely integrable system has a \emph{singularity} at
every point where this linear independence fails to hold.  It is
precisely at the singularities where the most interesting, and most
complicated, dynamical features of the system are displayed.  An
important class of completely integrable systems, with well\--behaved
singularities, are those given by Hamiltonian $n$\--torus actions on
symplectic $2n$\--manifolds. These actions have a momentum map with
$n$ components $f_1,\ldots,f_n$, and these components always form a
completely integrable system. A remarkable structure theory by Atiyah
\cite{atiyah}, Guillemin\--Sternberg \cite{gs} and Delzant
\cite{delzant} exists for these systems, which are usually referred to
as \emph{toric systems}.

The study of completely integrable Hamiltonian systems is a vast and
active research area.  Two motivations to study such systems come
from: (i) \emph{Kolmogorov\--Arnold\--Moser (KAM) theory}: since
integrable systems are ``solvable'' in a precise sense, one expects to
find valuable information about the behavior of dynamical systems that
are obtained by \emph{small perturbations} of them, and then the
powerful KAM theory comes into play (see de la Llave's article
\cite{llave} for a summary of the main ideas of KAM theory) to deal
with the properties of the perturbations (persistence of
quasi\--periodic motions); (ii) \emph{the theory of singularities of
  fibrations $(f_1,\ldots, f_n) \colon M \to \mathbb{R}$ by the
  Fomenko school} \cite{bolsinov-fomenko-book}: the Fomenko school has
developed powerful and far reaching methods to study the topology of
singularities of integrable systems.  It is interesting to notice that
there is a relation between these two motivations, which has been
explored recently by Dullin\--V\~u Ng\d oc and Nguy{\^e}n Ti{\^e}n
Zung \cite{san-dullin, san-dullin2, zung-kolmogorov, zungmore2}.

In the present article we give an overview of our perspective of the
current state of the art of the symplectic geometry of completely
integrable systems, with a particular emphasis on the the recent
developments on semitoric integrable systems in dimension four.
Before this, we briefly review several preceding fundamentals results
due to Arnold, Atiyah, Carath{\'e}odory, Darboux, Delzant,
Duistermaat, Dufour, Eliasson, Guillemin, Liouville, Mineur, Molino,
Sternberg, Toulet and Nguy{\^e}n Ti{\^e}n Zung, some of which are a
key ingredient in the symplectic theory of semitoric integrable
systems. This article does not intend to be comprehensive in any way,
but rather it is meant to be a fast overview of the current research
in the subject; we hope we will convey some of the developments which
we consider most representative.  Our point of view is that of local
phase\--space analysis; it advocates for the use of local normal
forms, and sheaf theoretic methods, to prove global results by gluing
local pieces.

Some of the current activity on integrable systems is concerned with a
question of high interest to applied and pure mathematicians and
physicists.  The question is whether one can reconstruct an integrable
system that one does not a priori know, from observing some of its
properties.  E.g.  Kac's famous question: can you hear the shape of a
drum?  Kac's question in the context of integrable systems can be
formulated in the following way: can a completely integrable system be
recovered from the semiclassical joint spectrum of the corresponding
\emph{quantized} integrable system?  In order to study this question
one must complement Fomenko's \emph{topological} theory with a
\emph{symplectic} theory, which allows one to \emph{quantize} the
integrable system. Quantization (the process of assigning Hilbert
spaces and operators to symplectic manifolds and smooth real\--valued
functions) is a motivation of the present article but it is not its
goal; instead we lay down the symplectic geometry needed to study the
quantization of certain completely integrable systems, and we will
stop there; occasionally in the paper, and particularly in the last
section, we will make some further comments on quantization.  For a
basic reference on the so called \emph{geometric quantization}, see
for example Kostant\--Pelayo \cite{kostantpelayo}.

The authors have recently given a global \emph{symplectic}
classification of integrable systems with two degrees of
freedom\footnote{The number of degrees of freedom is half the
  dimension of the symplectic manifold.} which has no hyperbolic
singularities and for which one component of the system is a
$2\pi$\--periodic Hamiltonian \cite{pelayovungoc1, pelayovungoc2};
these systems are called \emph{semitoric}.  We devote sections
\ref{sec:semitoric}, \ref{sec:inv}, \ref{sec:semitoric2} this paper to
explain this symplectic classification. This symplectic classification
of semitoric integrable system described in this paper prepares the
ground for answering Kac's question in the context of semitoric
completely integrable systems.

Semitoric systems form an important class of integrable systems,
commonly found in simple physical models; a semitoric system can be
viewed as a Hamiltonian system in the presence of circular
symmetry. Perhaps the simplest example of a non\--compact non\--toric
semitoric system is the coupled spin\--oscillator model
$S^2\times\R^2$ described in~\cite[Section~6.2]{vungoc}, where $S^2$ is
viewed as the unit sphere in $\R^3$ with coordinates $(x,\,y,\,z)$,
and the second factor $\R^2$ is equipped with coordinates $(u,\, v)$,
equipped with the Hamiltonians $ J := (u^2+v^2)/2 + z$ and $H :=
\frac{1}{2} \, (ux+vy)$.  Here $S^2$ and $\mathbb{R}^2$ are equipped
with the standard area forms, and $S^2 \times \mathbb{R}^2$ with the
product symplectic form. The authors have carried out the quantization
of this model in \cite{example}.

In the aforementioned papers we combine techniques from classical
differential geometry, semiclassical analysis, and Lie theory; these
works are representative of our core belief that one can make definite
progress in the program to understand the symplectic and spectral
theory of integrable systems by combining techniques and ideas from
these areas.  This symplectic work in turns generalizes the celebrated
theory of Hamiltonian Lie group actions by Atiyah, Benoist, Delzant,
Guillemin, Kirwan, Sternberg and others, to completely integrable
systems. It is also intimately connected with several previous works
\cite{DuPe, duispel3, duispel2, duispel4, kirwan, pelayo, vungoc, vungoc0}.

While major progress has been made in recent times by many authors,
the theory of integrable systems in symplectic geometry is far from
complete at the present time, even in the case of integrable systems
with two degrees of freedom. For example, it remains to understand the
symplectic theory of integrable systems on $4$\--manifolds when one
allows hyperbolic singularities to occur. The presence of hyperbolic
singularities has a global effect on the system which makes describing
a set of global invariants difficult\footnote{starting with what we
  will call the ``polygon invariant'' and which encodes in some
  precise sense the affine structure induced by the system.}.  We do
not know at this time if this is even a feasible problem or whether
one can expect to give a reasonable classification extending the case
where no hyperbolic singularities occur.

Moreover, the current theory allows us to understand semitoric systems
with controlled behavior at infinity; precisely this means that the
$2\pi$\--periodic Hamiltonian is a proper map (the paper in the works \cite{PeRaVu2010}
is expected to address this case). The general case is open,
however.

Shedding light into these two questions would bring us a step
closer to understanding the symplectic geometry of \emph{general}
completely integrable systems with two degrees of freedom in dimension
$4$, which we view as one of the major and longstanding unsolved
problems in geometry and dynamics (and to which many people have done
contributions, several of which are mentioned in the present paper).
In addition, answers to these questions constitute another required
step towards a quantum theory of integrable systems on symplectic
manifolds.

One can find integrable systems in different areas of mathematics and
physics. For example, in the context of algebraic geometry a semitoric
system naturally gives a toric fibration with singularities, and its
base space becomes endowed with a singular integral affine structure.
Remarkably, these singular affine structures are of key importance in
various parts of symplectic topology, mirror symmetry, and algebraic
geometry -- for example they play a central role in the work of
Kontsevich and Soibelman \cite{KS}, cf. Section \ref{rmk:gross} for
further details.  Interesting semitoric systems also appear as
relevant examples in the theory of symplectic quasi-states, see
Eliashberg\--Polterovich \cite[page 3]{eliashberg}. Many aspects of
the global theory of semitoric integrable systems may be understood in
terms of singular affine structures, but we do not know at this time
whether all of the invariants may be expressed in terms of singular
affine structures (if this were the case, likely involving some
asymptotic behavior).

For mathematicians semitoric systems are the next natural class of
integrable systems to consider after toric systems. Semitoric systems
exhibit a richer, less rigid behavior than toric systems.  The
mathematical theory of semitoric systems explained in the last few
sections of this paper was preceded by a number of interesting works
by physicists and chemists working on describing energy-momentum
spectra of systems in the context of quantum molecular spectroscopy
\cite{Fi2009, Sa1995, Ch1999, As2010}. Physicists and chemists were
the first to become interested in semitoric systems. Semitoric systems
appear naturally in the context of quantum chemistry. Many groups have
been working on this topic, to name a few: Mark Child's group in
Oxford (UK), Jonathan Tennyson's at University College London (UK),
Frank De Lucia's at Ohio State University (USA), Boris Zhilinskii's at
Dunkerque (France), and Marc Joyeux's at Grenoble (France).

These physicists and chemists have asked whether one can one give a
finite collection of invariants characterizing systems of this nature.
The theory of semitoric systems described in the present paper was
largely motivated by this question, and fits into the broader
realization in the physics and chemistry communities that symplectic
invariants play a leading role in understanding a number of global
questions in molecular spectroscopy --- hence any mathematical
discovery in this direction will be of interest outside of a pure
mathematical context, see Stewart \cite{St2004}.

Direct applications of integrable systems can also be found in the
theory of geometric phases, nonholonomic mechanics, rigid body
systems, fluid mechanics, elasticity theory and plasma physics, and
have been extensively carried out by many authors, including Marsden,
Ratiu and their collaborators. The semiclassical aspects of integrable
systems have been recently studied in the book \cite{san-book} and the
article \cite{sansurveypaper}. In the book \cite{CuBa1997} singular
Lagrangian fibrations are treated from the point of view of classical
mechanics. Finally, we would like to point out Bolsinov\--Oshemkov's
interesting review article \cite{BoOs2006}, where for instance one can
find very interesting information about hyperbolic singularities.

The structure of the paper is as follows. In sections
\ref{sec:dynamics}, \ref{sec:systems}, \ref{sec:localsystems} and
\ref{sec:semilocalsystems} we summarize some of the most important
known results at a local and semiglobal level for Hamiltonian systems,
and motivate their study by presenting some influential results from
the theory of Hamiltonian Lie group actions due to Atiyah,
Guillemin\--Sternberg and Delzant. In section \ref{sec:semitoric} we
introduce semitoric systems in dimension four, and explain their
convexity properties. In sections \ref{sec:inv} and
\ref{sec:semitoric2} we introduce symplectic invariants for these
systems and explain the recent global symplectic classification of
semitoric systems given by the authors. In Section
\ref{sec:semitoric3} we briefly discuss some open problems.
\\
\\
\\
\emph{Acknowledgments}.  We thank the anonymous referees for many
helpful comments and remarks.  AP was partly supported by NSF
Postdoctoral Fellowship DMS-0703601, NSF Grant DMS-0965738, a Leibniz
Fellowship from the Oberwolfach Foundation, and MSRI and IAS
memberships. He gratefully acknowledges the Mathematical Sciences
Research Institute in Berkeley for the hospitality during the academic
year 2009\--2010, and the Mathematisches Forschungsinstitut
Oberwolfach for the hospitality during part of the summer of 2010.
SVN was partly supported by an ANR ``Programme Blanc''.  Finally,
thanks are due to Michael VanValkenburgh for his careful reading of a
preliminary version.

\section{Symplectic Dynamics} \label{sec:dynamics}

The unifying topic this paper is \emph{symplectic geometry}, which is
the mathematical language to clearly and concisely formulate problems
in classical physics and to understand their quantum counterparts (see
Marsden\--Ratiu classical textbook \cite{marsdenratiu} for a treatment
of classical mechanical systems). In the sense of Weinstein's creed,
symplectic geometry is of interest as a series of remarkable
``transforms" which connect it with several areas of semiclassical
analysis, partial differential equations and low\--dimensional
topology.

One may argue that symplectic manifolds are not the most general, or
natural, setting for mechanics. In recent times some efforts have been
made to study Poisson structures, largely motivated by the study of
coadjoint orbits. However only very few general results on integrable
systems are known in the context of Poisson manifolds.

\subsection{Symplectic manifolds}

A \emph{symplectic form} on a vector space $V$ is a non\--degenerate,
antisymmetric, bilinear map $V \times V \to \R$.  A \emph{symplectic
  manifold} is a pair $(M,\omega)$ where $M$ is a smooth manifold and
$\omega$ is symplectic form on $M$, i.e. smooth collection of
symplectic forms $\omega_p$, one for each tangent space $\op{T}_pM$,
which is globally closed in the sense that the differential equation
$\op{d}\!\omega=0$ holds.

The simplest example of a symplectic manifold is probably a surface of
genus $g$ with an area form. An important example is $
\mathbb{R}^{2n}$ with the form $\sum_{i=1}^n\op{d}\!x_i \wedge \op{d}\!y_i, $
where $(x_1,\,y_1,\ldots,\,x_n,\,y_n)$ are the coordinates in
$\mathbb{R}^{2n}$.

Symplectic manifolds are always even dimensional, so for example $
S^3, \, S^1 $ cannot be symplectic. They are also orientable, where
the volume form is given by $\omega \wedge \ldots (n \,\,
\textup{times})\, \ldots \wedge \omega=\omega^n$, if $\dim M=2n$, so
for example the Klein bottle is not a symplectic manifold.  Moreover,
symplectic manifolds are topologically ``non\--trivial'' in the sense
that if $M$ is compact, then the even dimensional de Rham cohomology
groups of $M$ are not trivial because $[\omega^{k}] \in
\op{H}^{2k}_{\textup{dR}}(M)$ defines a non\--vanishing cohomology
class if $k\le n$, i.e.  the differential $2$\--form $\omega^k$ is
closed but not exact (the proof of this uses Stokes' theorem, and is
not completely immediate). Therefore, the spheres $ S^4,\,
S^6,\,S^8,\, \ldots, S^{2N},...  $ cannot be symplectic. Symplectic
manifolds were locally classified by Darboux and the end of the XIX
century. He proved the following remarkable theorem.

\begin{theorem}[Darboux  \cite{darboux}]
  Near each point in $(M,\, \omega)$ there exists coordinates $(
  x_{1},y_{1},\ldots,x_n,y_n)$ in which the symplectic form $\omega$
  has the form $\omega\,\,=\,\,\sum_{i=1}^n \,\,\op{d}\!x_i \wedge
  \op{d}\!y_i$.
\end{theorem}

It follows from Darboux's theorem that symplectic manifolds have no
local invariants other than the dimension. This is a fundamental
difference with Riemannian geometry, where the curvature is a local
invariant.

\subsection{Dynamics of vector fields and torus actions}

A smooth vector field $\mathcal{Y}$ on a symplectic manifold
$(M,\,\omega)$ is \emph{symplectic} if its flow preserves the
symplectic form $\omega$, and it is \emph{Hamiltonian} if the system
$$
\omega(\mathcal{Y},\, \cdot)=\op{d}\!H \,\,\,\,\,\,
(\textup{Hamilton's PDEs})
$$
has a smooth solution $H \colon M \to \mathbb{R}$. If so, we use the
notation $\mathcal{Y}:=\mathcal{H}_H$, and call $H$ call the
\emph{Hamiltonian} or \emph{Energy Function}.

\begin{figure}[htbp]
  \begin{center}
    \includegraphics[width=5.5cm]{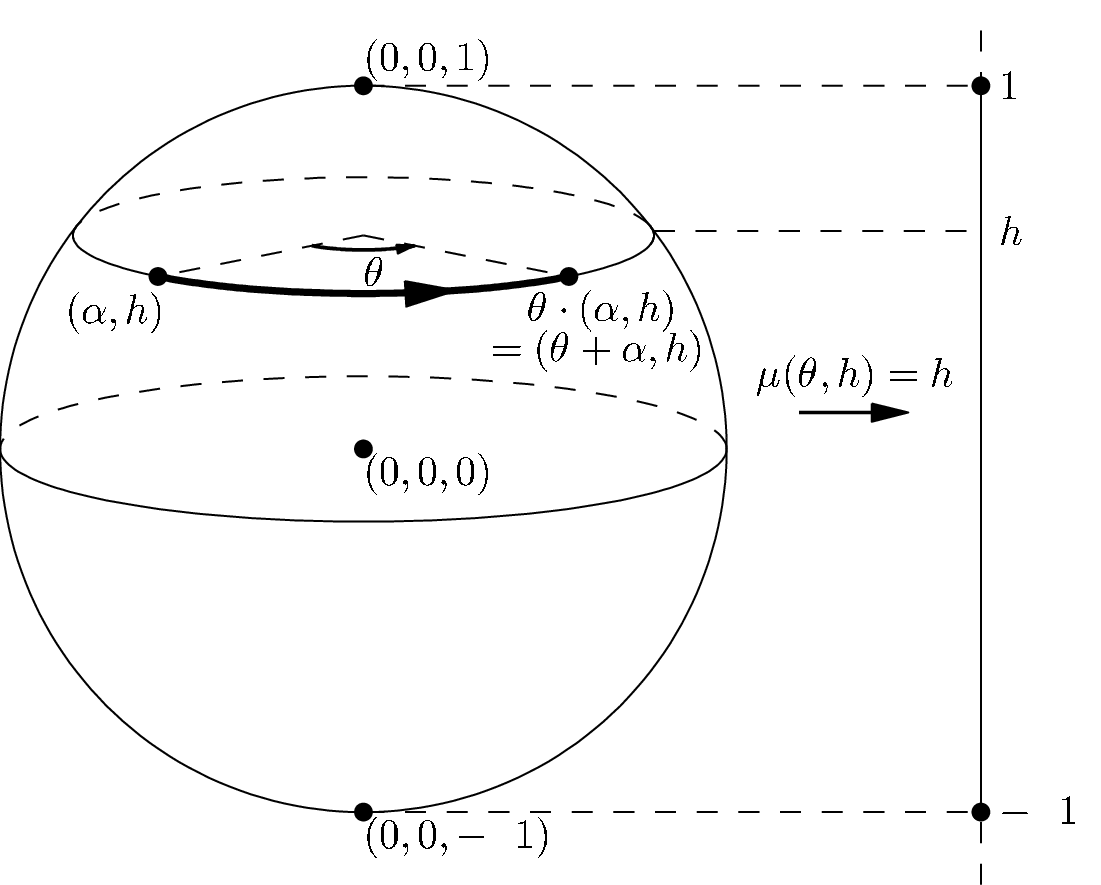}
    \caption{The momentum map for the $2$\--sphere $S^2$ is the height
      function $\mu(\theta,\, h)=h$. The image of $S^2$ under the
      momentum map $\mu$ is the closed interval $[-1,\,1]$. Note that
      as predicted by the Atiyah\--Guillemin\--Sternberg Theorem (see
      Theorem \ref{gs}), the interval $[-1,\,1]$ is equal to the image
      under $\mu$ of the set $\{(0,\,0,\,-1),\, (0,\, 0,\, 1)\}$ of
      fixed points of the Hamiltonian $S^1$\--action on $S^2$ by
      rotations about the vertical axis.}
    \label{fig1}
  \end{center}
\end{figure}

For instance, the vector field $\frac{\partial}{\partial \theta}$ on
$\mathbb{T}^2:=(S^1)^2$ is symplectic but not Hamiltonian ($\theta$ is
the coordinate on the first copy of $S^1$ in $\mathbb{T}^2$, for
example); on the other hand, the vector field
$\frac{\partial}{\partial \theta}$ on $S^2$ is Hamiltonian:
$\frac{\partial}{\partial \theta}=\mathcal{H}_{H}$ with
$H(\theta,\,h):=h$; here $(\theta,\, h)$ represents a point in the
unit sphere $S^2$ of height $h$ measured from the plane $z=0$ and
angle $\theta$ measured about the vertical axis, see Figure
\ref{fig1}.

Suppose that we have local Darboux coordinates
$(x_1,y_1,\ldots,x_n,y_n)$ near a point $m \in M$. Let $
\gamma(t):=(x_1(t),\,y_1(t),\ldots,x_n(t),\,y_n(t)) $ be an integral
curve of a smooth vector field $\mathcal{Y}$. Then
$\mathcal{Y}=\mathcal{H}_H$ for a local smooth function $H$ if and
only if
\[
\left\{ \begin{array}{rl}
    \frac{\op{d}\!y_i}{\op{d}\!t}(t)&\,\,=\,\,\,\,-\frac{\partial H}{\partial x_i}(\gamma(t)) \\
    \textup{\,} \\
    \frac{\op{d}\!x_i}{\op{d}\!t}(t)&\,\,=\,\,\,\,\,\,\,\,\,\,\,\frac{\partial
      H}{\partial y_i}(\gamma(t))
  \end{array}   
\right.
\]
It always holds that $H(\gamma(t))=\textup{const}$, i.e. that energy
is conserved by motion (Noether's Principle). If $\mathcal{Y}$ is
symplectic, these equations always have a local solution, but in order
for the vector field $\mathcal{Y}$ to be Hamiltonian (globally) one
must have that $\mathcal{Y}=\mathcal{H}_H$ on $M$, i.e. the function
$H$ has to be the same on each local Darboux chart. From a more
abstract point of view, this amounts to saying that the $1$\--form
$\omega(\mathcal{Y},\cdot)$ is always locally exact, but not
necessarily globally exact. So the obstruction to being exact lies in
the first de Rham cohomology group $\op{H}^1_{\textup{dR}}(M)$ of $M$;
if this group is trivial, then any smooth symplectic vector field on
$M$ is Hamiltonian.

Now suppose that we have a torus $T \simeq \mathbb{T}^k:=(S^1)^k$,
i.e.  a compact, connected abelian Lie group.  Let $X \in
\mathfrak{t}=\op{Lie}(T)$.  For $X$ in the Lie algebra $\mathfrak{t}$
of $T$ (i.e. we view $X$ as a tangent vector at the identity element
$1$ to $T$), there exists a unique \emph{homomorphism} $\alpha_X\colon
\R \to T$ with $\alpha_X(0)=1,\alpha_X'(0)=X$.  Define the so called
\emph{exponential map} $\op{exp} \colon \mathfrak{t} \to T$ by
$\op{exp}(X):=\alpha_X(1)$.  Using the exponential map, one can
generate many vector fields on a manifold from a given torus action
action. Indeed, for each $X \in \mathfrak{t}$, the vector field
$\mathcal{G}(X)$ on $M$ \emph{generated by $T$\--action} from $X$ is
defined by
\[
\mathcal{G}(X)_p :=\textcolor{black}{\text{tangent vector to}}
\underbrace{ t \mapsto
  \overbrace{\op{exp}(tX)}^{\textup{\textcolor{black}{curve in}}\,\,
    T}\cdot p}_{\textup{\textcolor{black}{curve in}}\,\, M\,\,
  \textcolor{black}{\textup{through}}\,\,p}
\textcolor{black}{\text{at}}\,\, t=0
\]
A $T$\--action on $(M,\, \omega)$ is \emph{symplectic} if all the
vector fields that it generates are \emph{symplectic}, i.e. their
flows preserve the symplectic form $\omega$. The $T$\--action is
\emph{Hamiltonian} if all the vector fields it generates are
\emph{Hamiltonian}, i.e.  they satisfy Hamilton's PDEs. Any symplectic
action on a simply connected manifold is Hamiltonian.

From a given Hamiltonian torus action, one can construct a special
type of map, which encodes information about the action -- this is the
famous \emph{momentum map}. The construction of the momentum map is
due to Kostant \cite{kostant1966} and Souriau \cite{souriau1966} (we
refer to Marsden\--Ratiu \cite[Pages 369, 370]{marsdenratiu} for the
history of the momentum map); the momentum map can be defined with
great generality for a Hamiltonian Lie group action.  The momentum map
was a key tool in Kostant's quantization lectures \cite{kostant1970}
and Souriau discussed it at length in his book
\cite{souriau1970}. Here we shall only deal with the momentum map in
the rather particular case of torus actions, in which the construction
is simpler.

Assume that $\op{dim}T=m$, $\op{dim}M=2n$.  Let $e_1,\ldots,e_m$ be
basis of the Lie algebra $\mathfrak{t}$. Let
$\mathcal{E}_1,\ldots,\mathcal{E}_m$ be the corresponding vector
fields. By definition of Hamiltonian action there exists a unique (up
to a constant) Hamiltonian $H_i$ such that $\omega(\mathcal{E}_i,
\cdot):=-\op{d}\!H_i$, i.e. $\mathcal{E}_i=\mathcal{H}_{H_i}$. Now we define the
\emph{momentum map} by
$$
\mu:=(H_1,\ldots,H_m) \colon M \to \mathbb{R}^m.
$$
The map $\mu$ is unique up to composition by an element of
$\op{GL}(m,\, \mathbb{Z})$ (because our construction depends on the
choice of a basis) and translations in $\mathbb{R}^m$ (because the
Hamiltonians are defined only up to a constant).

The simplest example of a Hamiltonian torus action is given by $S^2$ with the
rotational $S^1$\--action depicted in Figure \ref{fig1}. It is easy to
check that the momentum map for this action is the height function
$\mu(\theta, \,h)= h$.

In the US East Coast the momentum map has traditionally been called
``moment map'', while in the West Coast it has been traditional to use
the term ``momentum map''. In French they both reconcile into the term
``application moment''.

\subsection{Structure theorems for Hamiltonian actions}

Much of the authors' intuition on integrable systems was originally
guided by some remarkable results proven in the early 80s by Atiyah,
Guillemin\--Sternberg, and Delzant, in the context of Hamiltonian
torus actions.  The first of these results was the following
influential convexity theorem of Atiyah, Guillemin\--Sternberg.

\begin{theorem}[Atiyah \cite{atiyah}, Guillemin\--Sternberg
  \cite{gs}] \label{gs} If an $m$\--dimensional torus acts on a
  compact, connected $2n$\--dimensional symplectic manifold $(M,\,
  \omega)$ in a Hamiltonian fashion, the image $\mu(M)$ under the
  momentum map $ \mu:=(H_1,\ldots,\, H_m) \colon M \to \R^m $ is a
  convex polytope.
\end{theorem}

\begin{figure}[htbp]
  \begin{center}
    \includegraphics[width=5cm]{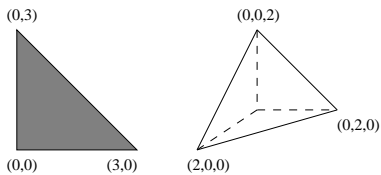}
    \caption{Delzant polytopes corresponding to the complex projective
      spaces $\mathbb{CP}^2$ and $\mathbb{CP}^3$ equipped with scalar
      multiples of the Fubini\--Study symplectic form.}
    \label{fig:AFF}
  \end{center}
\end{figure}

See figures \ref{fig1} and \ref{fig:AFF} for an illustration of the
theorem. Other remarkable convexity theorems were proven after the
theorem above by Kirwan \cite{kirwan} (in the case of compact,
non\--abelian group actions), Benoist \cite{benoist} (in the case when
the action is not necessarily Hamiltonian but it has some coisotropic
orbit) and Giacobbe \cite{giacobbe}. Convexity in the case of Poisson
actions has been studied by Alekseev, Flaschka\--Ratiu, Ortega\--Ratiu
and Weinstein \cite{Alek, FlaRat, OrtRat, Weinstein} among others.

Recall that a convex polytope in $\R^n$ is \emph{simple} if there are
$n$ edges meeting at each vertex, \emph{rational} if the edges meeting
at each vertex have rational slopes, i.e. they are of the form
$p+tu_i$, $0 \le t < \infty$, where $u_i \in \Z^n$, and \emph{smooth}
if the vectors $u_1,\ldots,u_n$ may be chosen to be a basis of $\Z^n$
(see Figure \ref{fig:AFF}). In the mid 1980s Delzant \cite{delzant}
showed the following classification result, which complements the
Atiyah\--Guillemin\--Sternberg convexity theorem.

\begin{theorem}[Delzant \cite{delzant}] \label{delzant} If an
  $n$\--dimensional torus acts effectively and Hamiltonianly on a
  compact, connected symplectic $2n$\--dimensional manifold $(M,\,
  \omega)$, the polytope in the Atiyah\--Guillemin\--Sternberg theorem
  is simple, rational and smooth, and it determines the symplectic
  isomorphism type of $M$, and moreover, $M$ is a toric variety in the
  sense of complex algebraic geometry. Starting from any simple,
  rational smooth polytope $\Delta \subset \R^m$ one can construct a
  compact connected symplectic manifold $(M_{\Delta},\,
  \omega_{\Delta})$ with an effective Hamiltonian action for which its
  associated polytope is $\Delta$.
\end{theorem}

By an isomorphism $\chi \colon (M_1,\, \omega_1) \to (M_2,\,
\omega_2)$ in Theorem \ref{delzant} we mean an equivariant
symplectomorphism such that $\chi^*\mu_2=\mu_1$, where $\mu_i$ is the
momentum map of $M_i$, $i=1,\,2$ (the map $\chi$ is an equivariant
symplectomorphism in the sense that it is a diffeomorphism which pulls
back the symplectic form $\omega_2$ to $\omega_1$ and commutes with
the torus actions). The manifolds in Delzant's theorem are called a
\emph{symplectic\--toric} manifolds or \emph{Delzant manifolds}. See
Duistermaat\--Pelayo \cite{duispel3} for a detailed study of the
relation between Delzant manifolds and toric varieties in algebraic
geometry.  In the context of symplectic geometry, motivated by
Delzant's results usually one refers to simple, rational smooth
polytopes as \emph{Delzant polytopes}.

Delzant's theorem tells us that from the point of view of symplectic
geometry, complex projective spaces endowed with the standard action
by rotations of a torus half the dimension of the corresponding
complex projective space are simple polytopes. More precisely,
consider the projective space $\mathbb{CP}^n$ equipped with a
$\lambda$ multiple of the Fubini--Study form
and the standard rotational action of $\mathbb{T}^n$ (for
$\mathbb{CP}^1=S^2$, we already drew the momentum map in Figure
\ref{fig1}). The complex projective space $\mathbb{CP}^n$ is a $2n$--dimensional
symplectic\--toric manifold, and one can check that the momentum map
is given by
$$
\mu(z)\,=(\frac{\lambda \,\, |z_1|^2}{\sum_{i=0}^n\,
  \,|z_i|^2},\ldots,\frac{\lambda \,\, |z_n|^2}{\sum_{i=0}^n\,
  \,|z_i|^2}).
$$
It follows that the momentum polytope equals the convex hull in
$\mathbb{R}^n$ of $0$ and the scaled canonical vectors $\lambda
e_1,\ldots,\lambda e_n$, see Figure \ref{fig:AFF}. Theorem \ref{delzant}
says that this polytope determines all the information about
$\mathbb{CP}^n$, the symplectic form and the torus action.

There have been many other contributions to the structure theory of
Hamiltonian torus actions, particularly worth noting is Karshon's
paper \cite{K} (see also \cite{K2}), where she gives a classification
of Hamiltonian circle actions on compact connected $4$\--dimensional
symplectic manifolds; we briefly review Karshon's result.  To a
compact, connected $4$\--dimensional symplectic manifold equipped with
an effective Hamiltonian $S^1$\--action (i.e. a so called
\emph{compact $4$\--dimensional Hamiltonian $S^1$\--space}), we may
associate a labelled graph as follows.  For each component $\Sigma$ of
the set of fixed points of the $S^1$\--action there is one vertex in
the graph, labelled by the real number $\mu(\Sigma)$, where $\mu
\colon M \to \mathbb{R}$ is the momentum map of the action. If
$\Sigma$ is a surface, then the corresponding vertex has two
additional labels, one being the symplectic area of $\Sigma$ and the
other one being the genus of $\Sigma$.

For every finite subgroup $F_k$ of $k$ elements of $S^1$ and for every
connected component $C$ of the set of points fixed by $F_k$ we have an
edge in the graph, labeled by the integer $k>1$. The component $C$ is
a $2$\--sphere, which we call a \emph{$F_k$\--sphere}. The quotient
circle $S^1/F_k$ rotates it while fixing two points, and the two
vertices in the graph corresponding to the two fixed points are
connected in the graph by the edge corresponding to $C$.

On the other hand, it was proven by Audin, Ahara and Hattori \cite{AH,
  A1, A2} that every compact $4$\--dimensional Hamiltonian
$S^1$\--space is isomorphic (meaning $S^1$\--equivariantly
diffeomorphic) to a complex surface with a holomorphic $S^1$\--action
which is obtained from $\mathbb{CP}^2$, a Hirzebruch surface or a
$\mathbb{CP}^1$\--bundle over a Riemann surface (with appropriate
circle actions), by a sequence of blow\--ups at the fixed points.

Let $A$ and $B$ be connected components of the set of fixed
points. The $S^1$\--action extends to a holomorphic action of the
group $\mathbb{C}^{\times}$ of non\--zero complex numbers. Consider
the time flow given by the action of subgroup $\op{exp}(t)$, $t \in
\mathbb{R}$.  We say that \emph{$A$ is greater than $B$} if there is
an orbit of the $\mathbb{C}^{\times}$\--action which at time
$t=\infty$ approaches a point in $A$ and at time $t=-\infty$
approaches a point in $B$.

Take any of the complex surfaces with $S^1$\--actions considered by
Audin, Ahara and Hattori, and assign a real parameter to every
connected component of the set of fixed points such that these
parameters are monotonic with respect to the partial ordering we have
just described. If the manifold contains two fixed surfaces then
assign a positive real number to each of them in such a way that the
difference between the numbers is given by a formula involving the
previously chosen parameters. Karshon proved \cite[Theorem 3]{K2} that
for every such a choice of parameters there exists an invariant
symplectic form and a momentum map on the complex surface such that
the values of the momentum map at the fixed points and the symplectic
areas of the fixed surfaces are equal to the chosen
parameters. Moreover, every two symplectic forms with this property
differ by an $S^1$\--equivariant diffeomorphism.  Karshon proved the
following classification result \emph{{\`a} la Delzant}.

\begin{theorem}[Karshon \cite{K}] \label{karshon:thm} If two compact
  Hamiltonian $S^1$\--spaces\footnote{I.e. two compact connected
    $4$\--dimensional manifolds equipped with an effective Hamiltonian
    $S^1$\--action.} have the same graph, then they are isomorphic
  (i.e. $S^1$\--equivariantly symplectomorphic). Moreover, every
  compact $4$\--dimensional Hamiltonian $S^1$\--space is isomorphic to
  one of the spaces listed in the paragraph above.
\end{theorem}

Again, in Theorem \ref{karshon:thm}, an isomorphism is an equivariant
symplectormphism which pulls back the momentum map on one manifold to
the momentum map on the other manifold. Theorem \ref{karshon:thm} has
useful consequences, for example: every compact Hamiltonian
$S^1$\--space admits a $S^1$\--invariant complex structure for which
the symplectic form is K{\"a}hler.

\subsection{Structure theorems for symplectic actions}

From the viewpoint of symplectic geometry, the situation described by
the momentum polytope is very rigid.  It is natural to wonder whether
the structure results Theorem \ref{gs} and Theorem \ref{delzant} for
Hamiltonian actions of tori persist in a more general context. There
are at least two natural ways to approach this question, which we
explain next.

\begin{figure}[htbp]
  \begin{center}
    \includegraphics[width=5.5cm]{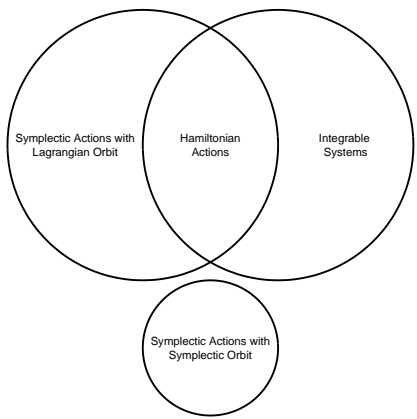}
    \caption{Hamiltonian torus actions may be viewed as a subclass of
      completely integrable Hamiltonian systems, which we study later
      in this paper, and a subclass of general symplectic actions
      which have some Lagrangian orbit. If the torus has dimension $n$
      and the manifold has dimension $2n$, Hamiltonian torus actions
      are a subclass of completely integrable systems.}
    \label{fig:fig1}
  \end{center}
\end{figure}

First, one can insist on having a compact group action, but not
requiring that the group acts in a Hamiltonian fashion. In other words
do the striking theorems above persist if the vector fields generated
by action have flows that preserve symplectic form i.e. are
symplectic, but Hamilton's PDEs have no solution i.e. the vector
fields are not Hamiltonian? Many easy examples fit this criterion,
e.g. take the $2$\--torus $\mathbb{T}^2$ with the standard area form
$\op{d}\! \theta \wedge \op{d}\! \alpha$ and with the
$\mathbb{S}^1$\--action on the $\theta$\--component; the basic vector
field $\frac{\partial}{\partial \theta}$ is symplectic but one can
easily check that it is non\--Hamiltonian.

Various works by Giacobbe \cite{giacobbe}, Benoist \cite{benoist},
Ortega\--Ratiu \cite{ortegaratiu}, Duistermaat\--Pelayo \cite{DuPe},
and Pelayo \cite{pelayo} follow this direction.  Benoist's paper gives
a convexity result for symplectic manifolds with coisotropic orbits,
Ortega\--Ratiu give a general symplectic local normal form theorem,
also studied by Benoist in the case that the orbits are coisotropic.
The papers by Duistermaat and Pelayo provide classifications
\emph{{\`a} la Delzant}. Let us briefly recall these classifications.

\begin{remark}
  Hamiltonian torus actions on compact manifolds always have fixed
  points (equivalently, the Hamiltonian vector fields generated by a
  Hamiltonian torus actions always have fixed points). Sometimes the
  condition of being Hamiltonian for vector fields can be detected
  from the existence of fixed points; this is in general a challenging
  question.

  The first result concerning the relationship between the existence
  of fixed points and the Hamiltonian character of vector fields
  generated by a $G$\--action is Frankel's celebrated theorem
  \cite{Frankel1959} which says that if the manifold is compact,
  connected, and K{\"a}hler, $G=S^1$, and the symplectic action has
  fixed points, then it must be Hamiltonian. Frankel's influential
  work has inspired subsequent research.  McDuff \cite[Proposition
  2]{McDuff1988} has shown that any symplectic circle action on a
  compact connected symplectic 4-manifold having fixed points is
  Hamiltonian. See Tolman\--Weitsman \cite[Theorem 1]{ToWe2000},
  Feldman \cite[Theorem 1]{Feldman2001}, \cite[Section 8]{KeRuTr2008},
  \cite{LuOp1995}, Giacobbe \cite[Theorem 3.13]{giacobbe},
  Duistermaat\--Pelayo \cite[Corollary 3.9]{DuPe}, Ginzburg
  \cite[Proposition 4.2]{Ginzburg1992}, Pelayo\--Tolman \cite{peltol}
  for additional results in the case of compact manifolds, and
  Pelayo\--Ratiu \cite{PelRat} for results in the case of
  non\--compact manifolds.
\end{remark}

Our next goal is to present a classification \emph{{\`a} la Delzant}
of symplectic torus actions that have some Lagrangian orbit; this in
particular includes all symplectic toric manifolds, because the
maximal (in the sense of dimension) orbits of a symplectic\--toric
manifold are Lagrangian.  An $n$\--dimensional submanifold $L$ of a
symplectic $2n$\--manifold $(M,\omega)$ is \emph{Lagrangian} if the
symplectic form $\omega$ vanishes on $L$. For example, the orbits of
the $S^1$\--action by rotations on $S^2$ in Figure \ref{fig1} are
Lagrangian because an orbit is given by $h=\textup{const}$ for some
constant and the symplectic form is $\op{d}\!\theta \wedge \op{d}\!
h$, which clearly vanishes when $h$ is constant.  So the maximal
orbits of the standard symplectic $2$\--sphere are Lagrangian.

A famous example of a symplectic manifold with a $2$\--torus action
for which all the orbits are Lagrangian is the the
\emph{Kodaira\--Thurston manifold}. It is constructed as follows: let
$(j_1, \,j_2) \in \mathbb{Z}^2$ act on $\mathbb{R}^2$ by the inclusion
map  (i.e. $(j_1,\,j_2) \cdot (x_1,\,y_1)=(j_1+x_1,\,
j_2+y_1)$), on $\mathbb{T}^2$ by the $2$ by $2$ matrix with entries $a_{11}=a_{22}=1$,
$a_{12}=j_2$, $a_{21}=0$, and on the product $\mathbb{R}^2 \times
\mathbb{T}^2$ by the diagonal action. This diagonal action gives rise
to a torus bundle over a torus $\mathbb{R}^2 \times_{\mathbb{Z}^2}
\mathbb{T}^2$, the total space of which is compact and connected. The
product symplectic form $ \textup{d}x_1 \wedge \textup{d}y_1 +
\textup{d}x_2 \wedge \textup{d}y_2$ on $\mathbb{R}^2 \times
\mathbb{Z}^2$ descends to a symplectic form on $\mathbb{R}^2
\times_{\mathbb{Z}^2} \mathbb{T}^2$.

Moreover, one can check that $\mathbb{T}^2$ acts symplectically on
$\mathbb{R}^2 \times_{\mathbb{Z}^2} \mathbb{T}^2$, where the first
circle of $\mathbb{T}^2$ acts on the left most component of
$\mathbb{R}^2$, and the second circle acts on the right most component
of $\mathbb{T}^2$ (one can check that this is indeed a well\--defined,
free symplectic action).  Because the action is free, it does not have
fixed points, and hence it is not Hamiltonian (it follows from the
Atiyah\--Guillemin\--Sternberg theorem that Hamiltonian actions always
have some fixed point).  All the orbits of this action are Lagrangian
submanifolds, because both factors of the symplectic form vanish since
each factor has a component which is zero because it is the
differential of a constant.

\begin{figure}[htbp]
  \begin{center}
    \includegraphics[height=7.5cm, width=7cm]{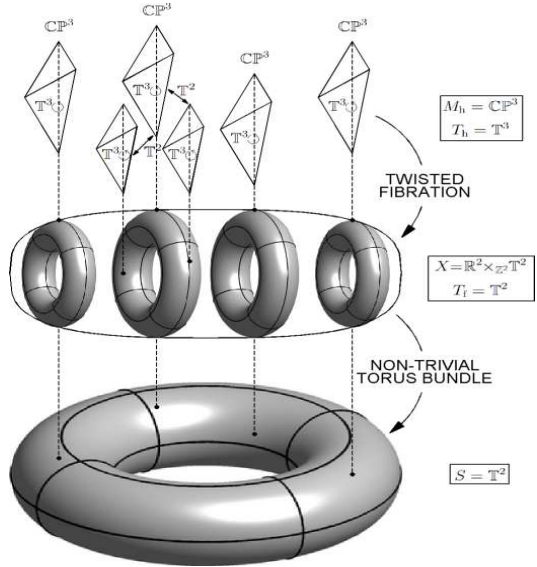}
    \caption{$10$\--dimensional symplectic manifolds with a torus
      action with Lagrangian orbits.  Vector field generated by
      $T$\--action is ``twist'' of
      $\mathcal{Y}=(\mathcal{Y}_{\textup{h}},\mathcal{Y}_{\textup{f}})$,
      where $\mathcal{Y}_{\textup{h}}$ Hamiltonian on $\mathbb{CP}^3$,
      $\mathcal{Y}_{\textup{f}}$ symplectic on
      $\mathbb{R}^2\times_{\mathbb{Z}^2}\mathbb{T}^2$.}
    \label{fig:???}
  \end{center}
\end{figure}

Another example is $\mathbb{T}^2 \times S^2$ equipped with the form
$\textup{d}x \wedge \textup{d}y +\textup{d}\theta \wedge \textup{d}h$,
on which the $2$\--torus $\mathbb{T}^2$ acts symplectically, one
circle on each factor.  This action has no fixed points, so it is not
Hamiltonian.  It is also not free. The free orbits are Lagrangian.
All of these examples, fit in the following theorem.

\begin{theorem}[Duistermaat-Pelayo \cite{DuPe}] \label{dp} Assume that
  a torus $\mathbb{T}^m$ of dimension $m$ acts effectively and
  symplectically on a compact, connected symplectic $2m$\--manifold
  $(M,\omega)$ with some Lagrangian orbit. Then $\mathbb{T}^m$
  decomposes as a product of two subtori
  $\mathbb{T}^m=T_{\textup{h}}T_{\textup{f}}$, where $T_{\textup{h}}$
  acts Hamiltonianly on $M$ and $T_{\textup{f}}$ acts freely on $M$,
  and there is two\--step fibration $M \to X \to S$, where $M$ is the
  total space of a fibration over $X$ with fibers symplectic\--toric
  manifolds $(M_{\textup{h}},T_{\textup{h}})$, and $X$ is a
  $T_{\textup{f}}$\--bundle over a torus $S$ of dimension
  $m-\textup{dim}\,T_{\textup{h}}$.
\end{theorem}

In this theorem $X$ is a symplectic homogenous space for the twisted
group $T \times \mathfrak{t}^*$.  The formulation of this theorem in
\cite{DuPe} is completely explicit, but it is too involved to be
described here. In particular, the formulation contains a complete
symplectic classification in terms of six symplectic invariants
(e.g. the Chern class of the fibration, the Hamiltonian torus
$T_{\textup{h}}$, the polytope corresponding to the Hamiltonian action
of $T_{\textup{h}}$ etc).  This classification includes Delzant's
classification (stated previously in the paper as Theorem
\ref{delzant}), which corresponds to the case of
$T_{\op{h}}=\mathbb{T}^m$ and $T_{\op{f}}$ trivial; in this case five
of the invariants do no appear, the only invariant is the
polytope. Note that the ``opposite'' situation occurs when
$T_{\op{f}}=\mathbb{T}^m$ and $T_{\op{h}}$ trivial (eg. the
Kodaira\--Thurston manifold), and in this case Theorem \ref{dp} says
that $M$ is a torus bundle over a torus with Lagrangian fibers.  An
example of a manifold which fits in Theorem \ref{dp} is the family of
10-dimensional twisted examples with Lagrangian orbits illustrated in
Figure \ref{fig:???}.

A classification theorem in the case when there exists a maximal
symplectic orbit (i.e.  an orbit on which the symplectic form
restricts to a symplectic form) was proven in \cite{pelayo}. In this
same paper, a classification of symplectic actions of $2$\--tori on
compact connected symplectic $4$\--manifolds was given building on
this result and Theorem \ref{dp}, but the description is involved for
the purposes of the present article.

\begin{theorem}[Pelayo \cite{pelayo}, Duistermaat\--Pelayo
  \cite{duispel2}] \label{DuPe} A compact connected symplectic
  $4$\--manifold $(M, \omega)$ equipped with an effective symplectic
  action of a $2$\--torus is isomorphic (i.e. equivariantly
  symplectomorphic) to one, and only one, of the spaces in the table:
  \begin{center}
    \begin{tabular}{ |clc|clc|clc|}
      \hline 
      {\bf {\tiny SPACE}}   & {\bf {\tiny ACTION}}   &  
      {\bf {\tiny MAXIMAL ORBITS}}      & {\bf {\tiny HAMILTONIAN?}} & 
      {\bf {\tiny INVARIANT COMPLEX?}} & {\bf {\tiny K{\"A}HLER?} }\\    \hline
      \textcolor{black}{{Toric}} & Fixed points & Lagrangian & Yes                     &  Yes &  Yes                                   \\ \hline
      $M \to \Sigma $ & Locally Free & Symplectic & No & Yes & Yes      \\ \hline
      $M \to \mathbb{T}^2$& Free & Lagrangian & No & Yes & No                                           \\   \hline
      $\mathbb{T}^2 \times S^2$ & Else & Lagrangian & No & Yes & Yes \\    \hline    
    \end{tabular}
  \end{center}
  The first item is a symplectic\--toric manifold with its standard
  Hamiltonian $2$\--torus action, the second item is an orbifold
  $2$\--torus bundle over a $2$\--dimensional compact connected
  orbifold, the third item is a $2$\--torus bundle over a $2$\--torus.
\end{theorem}

The first four columns in Theorem \ref{DuPe} were proven in
\cite{pelayo}, and the last two columns were proven in
Duistermaat\--Pelayo \cite{DuPe}; the article \cite{DuPe} is based on
Kodaira's seminal work \cite[Theorem 19]{kodaira} of 1961 on complex
analytic surfaces. Moreover, the two middle items above are completely
explicit and classified in terms of five symplectic invariants,
cf. \cite[Theorem 8.2.1]{pelayo}.

\begin{example}
Let us spell out the space on the third row in the table given in Theorem \ref{DuPe} 
more concretely, and we
refer to  \cite[Section 8]{pelayo} for the construction of the second row. The construction
which we present next is self\--contained and provides a source of 
many inequivalent examples. Let $T$ be a 
$2$\--dimensional torus. 
Let $T_{{\Z}}$ be the kernel of the exponential 
mapping $\textup{exp}:\mathfrak{t}\to T$.
\begin{itemize}
\item[a)]
For any choice of
\begin{itemize}
\item[i)]
a discrete cocompact subgroup $P$ of $\mathfrak{t}^*$, and
\item[ii)]
a non\--zero antisymmetric bilinear mapping $c \colon\mathfrak{t}^* \times \mathfrak{t}^* \to \mathfrak{t}$
such that $c(P\times P) \subset T_{\Z}$, 
\end{itemize}
let 
$\iota \colon P \to T \times \mathfrak{t}^*$ be
given by 
$
\zeta=\zeta_1\epsilon_1+\zeta_2\epsilon_2 \mapsto (
\op{e}^{-1/2 \,\zeta_1\zeta_2 \,c(\epsilon_1,\, \epsilon_2)} ,\,\zeta),
$ 
where $\epsilon_1,\,\epsilon_2$ is a $\Z$\--basis of $P$.
The mapping $\iota$ is a homomorphism onto a discrete cocompact subgroup of $T\times \got{t}^*$
with respect to the non\--standard standard group structure given by 
$$
(t,\,\zeta )\, (t',\,\zeta ')=(t\, t'\,\op{e}^{-c(\zeta ,\,\zeta ')/2},\,\zeta +\zeta ').
$$
Equip $T \times \mathfrak{t}^*$ with the standard cotangent bundle symplectic form.
Then 
$(T \times \mathfrak{t}^*)/\iota(P)$ equipped with the action of $T$
which comes from the action of $T$ by translations on the left factor of $T \times \mathfrak{t}^*$,
and where the symplectic form on $(T \times \mathfrak{t}^*)/\iota(P)$ is the $T$\--invariant form
induced by the symplectic form on $T \times \mathfrak{t}^*$, 
is a compact, connected symplectic $4$\--manifold on which $T$ acts freely and for which the
$T$\--orbits are Lagrangian $2$\--tori.  Theorem \ref{DuPe} implies that two symplectic manifolds 
constructed in this way
are isomorphic (i.e. equivariantly symplectomorphic) if and only if the corresponding
cocompact subgroups and the corresponding bilinlear forms are equal.
\item[b)]
For any choice of 
\begin{itemize}
\item[i)]
a discrete cocompact subgroup $P$ of $\mathfrak{t}^*$, and
\item[ii)]
a homomorphism $\tau \colon P \to T$, $\zeta \mapsto \tau_{\zeta}$,
\end{itemize}
let $\iota \colon P \to T \times \mathfrak{t}^*$ be given by $\zeta \mapsto (\tau_{\zeta}^{-1},\,\zeta)$. 
The mapping $\iota$ is a homomorphism onto a discrete cocompact subgroup of $T\times \got{t}^*$
with respect to the standard group structure.
Equip $T \times \mathfrak{t}^*$ with the standard cotangent bundle symplectic form.
Then 
$(T \times \mathfrak{t}^*)/\iota(P)$ equipped with the action of $T$
which comes from the action of $T$ by translations on the left factor of $T \times \mathfrak{t}^*$,
and where the symplectic form on $(T \times \mathfrak{t}^*)/\iota(P)$ is the $T$\--invariant form
induced by the symplectic form on $T \times \mathfrak{t}^*$, 
is a compact, connected symplectic $4$\--manifold on which $T$ acts freely 
with
$T$\--orbits Lagrangian $2$\--tori.  Theorem \ref{DuPe} implies that two symplectic manifolds constructed in
this way are $T$\--equivariantly
symplectomorphic if and only if the corresponding cocompact groups $P$ 
and the corresponding equivalence classes $\tau \cdot {\textup{exp}(\op{Sym}|_P)} \in \mathcal{T}$ are equal. 
Here $\op{exp} \colon \op{Hom}(P,\,\got{t}) \to
\op{Hom}(P,\,T)$ is the exponential map of the Lie group $\op{Hom}(P,\,T)$
and $\op{Sym}|_P \subset \op{Hom}(P,\,\got{t})$ is the space of restrictions
$\alpha|_P$ of linear maps $\alpha \colon \got{t}^* \to \got{t}$, $\xi \mapsto \alpha_{\xi}$,
which are symmetric in the sense that for all $\xi,\, \xi' \in {\got t}^*$,
$\xi(\alpha_{\xi'})-\xi'(\alpha_{\xi}) = 0$.
\end{itemize}
In both cases above the projection mapping $(T \times \got{t}^*)/\iota(P) \to \got{t}^*/P$
is a principal $T$\--bundle over the torus $\got{t}^*/P$ with Lagrangian fibers (the $T$\--orbits).
These spaces $(T \times \got{t}^*)/\iota(P)$ are all the possible cases that can occur as the third item in Theorem \ref{DuPe}.
\end{example}

Another natural generalization of a Hamiltonian torus action is the
notion of a completely integrable system, or more generally, of a
Hamiltonian system.  Probably the most fundamental difference between
the theory of Hamiltonian torus actions on compact manifolds and the
theory of completely integrable Hamiltonian systems can be seen
already at a local level.  Completely integrable systems have in
general singularities that are quite difficult to understand from a
topological, dynamical and analytic view\--point.  The singularities
of Hamiltonian torus actions occur at the lower dimensional orbits
only, and are tori of varying dimensions, but case of integrable
systems will usually exhibit a wider range of singularities such as
pinched tori, eg see Figure \ref{figure2}.  Indeed, it is only
recently that some of these singularities are beginning to be
understood in low dimensions. The rest of this paper is focused on the
local and global aspects of the symplectic geometry of completely
integrable systems.

\section{Completely Integrable Systems} \label{sec:systems}

\subsection{Hamiltonian integrable systems} \label{se}
Let $(M, \, \omega)$ be a $2n$\--dimensional symplectic manifold.  The
pair consisting of the smooth manifold and a classical observable
$H\colon M \to \R$ in $\op{C}^{\infty}(M)$ is called a
\emph{Hamiltonian system}.

A famous example of a Hamiltonian system is the \emph{spherical
  pendulum}, which is mathematically described as the symplectic
cotangent bundle $(\op{T}^*S^2, \omega_{\op{T}^*S^2})$ of the unit
sphere equipped with the Hamiltonian
\[
H(\underbrace{\theta,\phy}_{\tiny \textup{sphere}},
\underbrace{\xi_\theta,\xi_\phy}_{\textup{fiber}})=
\underbrace{\frac{1}{2}\left(\xi_\theta^2 +
    \frac{1}{\sin^2\theta}\xi_\phy^2\right)}_{\textup{kinetic
    energy}}+ \underbrace{\cos\theta}_{\textup{potential}}.
\]
Here $(\theta,\phy)$ are the standard spherical angles ($\phy$ stands
for the rotation angle around the vertical axis, while $\theta$
measures the angle from the north pole), and $(\xi_\theta,\xi_\phy)$
are the cotangent conjugate variables. The function $H$ is smooth on
$T^*S^2$ (the apparent singularity $1/\sin^2\theta$ is an artifact of
the spherical coordinates).

A classical observable $H$ gives rise to the \emph{Hamiltonian vector
  field} $\mathcal{H}_H$ on $M$ defined uniquely by $
\omega(\mathcal{H}_H,\, \cdot)=\op{d}\!H.  $ The algebra
$\op{C}^{\infty}(M)$ of classical observables, which we have been
calling \emph{Hamiltonians}, comes naturally endowed with the Poisson
bracket: $\{J,\,H\}:=\omega(\mathcal{H}_{J},\, \mathcal{H}_{H}).$ As a
derivation, $\mathcal{H}_H$ is just the Poisson bracket by $H$; in
other words the evolution of a function $f$ under the flow of
$\mathcal{H}_H$ is given by the equation $ \dot{f}=\{H,\,f\}.  $

An \emph{integral} of the Hamiltonian $H$ is a function which is
invariant under the flow of $\mathcal{H}_H$, i.e. a function $f$ such
that $\{H,\,f\}=0$.  The Hamiltonian $H$ is said to be
\emph{completely integrable} if there exists $n-1$ independent
functions $f_2,\dots,f_n$ (\emph{independent} in the sense that the
differentials $\op{d}_m\!H, \op{d}_m\!f_2,\,\dots,\,\op{d}_m\!f_n$ are
linearly independent at almost every point $m \in M$) that are
integrals of $H$ and that pairwise Poisson commute, i.e.  $\{H,\,
f_i\}=0$ and $\{f_i,\, f_j\}=0.$ For example, the spherical pendulum
$(\op{T}^*\!S^2, \omega_{\op{T}^*\!S^2},H)$ is integrable by
considering the vertical angular momentum, which is the function
$f_2(\theta,\varphi,\xi_\theta,\xi_\phy)= \xi_\phy$.  In the abstract
definition of a completely integrable system, it is clear that $H$
does not play a distinguished role among the functions
$f_2,\dots,f_n$.  The point of view in this paper will always be to
consider, as a whole, a collection of such functions. The integer $n$
is traditionally called the number of degrees of freedom of the system.

\begin{definition}
  A \emph{completely integrable system} on the $2n$\--dimensional
  symplectic manifold $M$, compact or not, is a collection of $n$
  Poisson commuting functions $f_1,\dots,f_n \in \op{C}^{\infty}(M)$
  which are independent.
\end{definition}

An important class of completely integrable systems are those given by
Hamiltonian $n$\--torus actions on symplectic $2n$\--manifolds
(i.e. symplectic\--toric manifolds).  These actions have a momentum
map with $n$ components $f_1,\ldots,f_n$, and these components always
form a completely integrable system in the sense of the definition
above.

\subsection{Singularities and regular points}
From a topological, analytical and dynamical view\--point, the most
interesting features of the completely integrable system on a
symplectic manifold are encoded in the ``singular'' fibers of the
momentum map $F=(f_1,\,\dots,\,f_n):M\to \R^n$, and in their
surrounding neighborhoods.

\begin{figure}[htbp]
  \begin{center}
    \includegraphics[width=14cm]{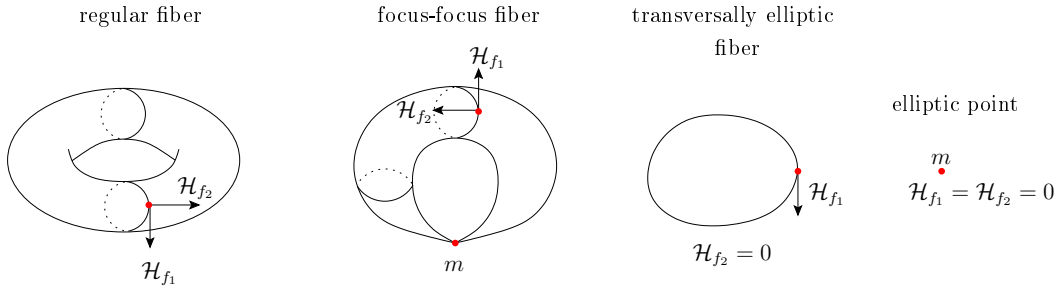}
    \caption{The figures show some possible singularities of a
      completely integrable system. On the left most figure $m$ is a
      regular point (rank 2); on the second figure $m$ is a
      focus\--focus point (rank $0$); on the third one $m$ is a
      transversally elliptic singularity (rank $1$); on the right most
      figure $m$ is a elliptic elliptic point.}
    \label{figure2}
  \end{center}
\end{figure}

A point $m\in M$ is called a \emph{regular point} if $\op{d}_mF$ has
rank $n$.  A point $c \in \R^n$ is a \emph{regular value} if the fiber
$F^{-1}(c)$ contains only regular points.  If $c$ is a regular value,
the fiber $F^{-1}(c)$ is called a \emph{regular fiber}.  A point $m
\in M$ is a \emph{critical point}, or a \emph{singularity}, if
$\op{d}_mF$ has rank strictly less than $n$. Geometrically this means
that the Hamiltonian vector fields generated by the components of $F$
are linearly dependent at $m$, see Figure \ref{figure2}. A fiber
$F^{-1}(c)$ is a \emph{singular fiber} if it contains at least one
critical point, see figure \ref{fig:pinched}.

It follows from the definition of a completely integrable system
(simply follow the flows of the Hamiltonian vector fields) that if $X$
is a connected component of a regular fiber $F^{-1}(c)$, where we
write
$$
F:=(f_1,\ldots,f_n) \colon M \to \mathbb{R}^n,
$$
and if the vector fields $\mathcal{H}_{f_1},\ldots, \mathcal{H}_{f_n}$
are complete on $F^{-1}(c)$, then $X$ is diffeomorphic to
$\R^{n-k}\times \T^k$.  Moreover, if the regular fiber ${F}^{-1}(c)$
is compact, then $\mathcal{H}_{f_1}, \ldots, \mathcal{H}_{f_n}$ are
complete, and thus the component $X$ is diffeomorphic to
$\mathbb{T}^n$; this is always true of for example some component
$f_i$ is proper.

The study of singularities of integrable systems is fundamental for
various reasons. On the one hand, because of the way an integrable
system is defined in terms $n$ smooth functions on a manifold, it is
expected (apart from exceptional cases) that singularities will
necessarily occur. On the other hand these functions define a
dynamical system such that their singularities correspond to fixed
points and relative equilibria of the system, which are of course one
of the main characteristics of the dynamics.

As a remark for those interested in semiclassical analysis we note
that from a semiclassical viewpoint, we know furthermore that
important wave functions such as eigenfunctions of the quantized
system have a microsupport which is invariant under the classical
dynamics; therefore, in a sense that we shall not present here (one
should talk about semiclassical measures), they concentrate near
certain singularities\footnote{those called hyperbolic.} (see for
instance~\cite{colin-p} and the work of
Toth~\cite{toth-expected}). This concentration entails not only the
growth in norm of eigenfunctions (see for
instance~\cite{toth-zelditch-1}) but also a higher local density of
eigenvalues~\cite{colin-p2,san-focus,san-colin}).

Let $f_1,\,\ldots,\,f_n$ define an integrable system on a symplectic
manifold $M$, and let $F$ be the associated momentum map. Suppose that
$F$ is a proper map so that the regular fibers of $F$ are
$n$\--dimensional tori.  Indeed, Liouville proved in 1855
\cite{liouville} that, locally, the equations of motion defined by any
of the functions $f_i$ are integrable by quadratures. This holds in a
neighborhood of any point where the differentials $\op{d}\!f_j$ are
linearly independent.

A pleasant formulation of Liouville's result, due to Darboux and
Carath{\'e}odory, says that there exist canonical coordinates $(x,\,
\xi)$ in which the functions $f_j$ are merely the ``momentum
coordinates'' $\xi_j$.  In 1935 Henri Mineur \cite{mineur}
stated\footnote{J.J. Duistermaat has pointed out \cite{DUpersonal}
  some gaps in Mineur's proof.}  in the special case of $\R^n \times
\R^n$ that if $\Lambda$ is a compact level set of the momentum map $F
= (f_1, \ldots, f_n)$, then $\Lambda$ is a torus. Moreover, there
exist symplectic coordinates $(x, \, \xi)$, where $x$ varies in the
torus $\T^n = \R^n/\Z^n$ and $\xi$ varies in a neighborhood of the
origin in $\R^n$, in which the functions $f_j$ depend only on the
$\xi$\--variables. In geometric terms, the system is symplectically
equivalent to a neighborhood of the zero section of the cotangent
bundle $\op{T}^*(\mathbb{T}^n)$ equipped with the integrable system
$(\xi_1,\ldots,\xi_n)$.  This result, proved in the general case in
1963 by Arnold \cite{arnold} (Arnold was not aware of Mineur's work),
is known as the \emph{action\--angle theorem} or the
\emph{Liouville\--Arnold theorem}. The tori are the famous
\emph{Liouville tori}. Although Liouiville's theorem has been
originally attached to this theorem for some time, we are not aware of
Liouville having contributed to this result; we thank J.J. Duistermaat
for pointing this out to us \cite{DUpersonal}.

We will study these and other results in more detail in the following
two sections.

\section{Local theory of completely integrable
  systems} \label{sec:localsystems}

\subsection{Local model at regular points} Let $(f_1,\dots,f_n)$ be a
completely integrable system on a $2n$-symplectic manifold $M$, with
momentum map $F$.  By the local submersion theorem, the fibers
$F^{-1}(c)$ for $c$ close to $F(m)$ are locally $n$-dimensional
submanifolds near a regular point $m$. The local structure of regular
points of completely integrable systems is simple:

\begin{theorem}[Darboux-Carath{\'e}odory]
  Let $(f_1,\dots,f_n)$ be a completely integrable system on a
  $2n$\--symplectic manifold $M$, with momentum map $F$.  If $m$ is
  regular, $F$ is symplectically conjugate near $m$ to the linear
  fibration $(\xi_1,\dots,\xi_n)$ on the symplectic space $\RM^{2n}$
  with coordinates $(x_1,\dots,x_n,\xi_1,\dots,\xi_n)$ and symplectic
  form $\sum_i \op{d}\!\xi_i\wedge \op{d}\!x_i$.
\end{theorem}

In other words, the Darboux-Carath{\'e}odory theorem says that there
exists smooth functions $\phi_1,\dots,\phi_n$ on $M$ such that $
(\phi_1, \dots, \phi_n, f_1, \dots, f_n) $ is a system of canonical
coordinates in a neighborhood of $m$.  In principle the name of
Liouville should be associated with this theorem, since well before
Darboux and Carath{\'e}odory, Liouville gave a nice explicit formula
for the functions $\phi_j$. This result published in
1855~\cite{liouville} explains the local integration of the flow of
any completely integrable Hamiltonian (possibly depending on time)
near a regular point of the foliation in terms of the \emph{Louville
  1-form} $\sum_i\xi_i\op{d}\!x_i$. In this respect it implies the
Darboux-Carath{\'e}odory theorem, even if Liouville's formulation is
more complicated.

\subsection{Local models at singular points}

One can approach the study of the singularities of Hamiltonian systems
in two different ways: one can analize the flow of the vector fields
--- this is the ``dynamical systems'' viewpoint --- or one can study
of the Hamiltonian functions themselves --- this is the ``foliation''
perspective.

In the case of\emph{ completely} integrable systems, the dynamical and
foliation points of view are equivalent because the vector fields of
the $n$ functions $f_1,\dots,f_n$ form a basis of the tangent spaces
of the leaves of the foliation $f_i=\textup{const}_i$, at least for
regular points.  The foliation perspective usually displays better the
geometry of the problem, and we will frequently use this
view\--point. However, the foliations we are interested in are
\emph{singular}, and the notion of a singular foliation is already
delicate. Generally speaking these foliations are of
Stefan-S{\"u}{\ss}mann type~\cite{stefan}~: the leaves are defined by
an integrable distribution of vector fields. But they are more than
that: they are Hamiltonian, and they are \emph{almost regular} in the
sense that the singular leaves cannot fill up a domain of positive
measure.

\subsubsection{Non\--degenerate critical points} \label{subsec:will}

In singularity theory for differentiable functions, ``generic''
singularities are Morse singularities. In the theory of completely
integrable systems there exists a natural analogue of the notion of
Morse singularities (or more generally of Morse-Bott singularities if
one allows critical submanifolds). These so-called
\emph{non-degenerate} singularities are now well defined and
exemplified in the literature, so we will only recall briefly the
definition in Vey's paper \cite{vey}.

Let $F=(f_1,\dots,f_n)$ be a completely integrable system on $M$.  A
fixed point $m\in M$ is called \emph{non-degenerate} if the Hessians
$\op{d}^2_mf_j$ span a Cartan subalgebra of the Lie algebra of
quadratic forms on the tangent space $\op{T}_m M$ equipped with the
linearized Poisson bracket. \label{def:nondegpoi}
This definition applies to a fixed point; more generally: if
$\op{d}_mF$ has corank $r$ one can assume that the differentials $
\op{d}_mf_1,\dots,\op{d}_mf_{n-r} $ are linearly independent; then we
consider the restriction of $f_{n-r+1},\dots,f_n$ to the symplectic
manifold $\Sigma$ obtained by local symplectic reduction under the
action of $f_1,\dots,f_{n-r}$. We shall say that $m$ is
\emph{non-degenerate} (or \emph{transversally non-degenerate})
whenever $m$ is a non-degenerate fixed point for this restriction of
the system to $\Sigma$.

In order to understand the following theorem one has to know the
linear classification of Cartan subalgebras of
$\textup{sp}(2n,\RM)$. This follows from the work of
Williamson~\cite{williamson}, which shows that any such Cartan
subalgebra has a basis build with three type of blocks: two
uni-dimensional ones (the elliptic block: $q=x^2+\xi^2$ and the real
hyperbolic one: $q=x\xi$) and a two-dimensional block called
focus-focus: $ \,\,\,\,\, q_1=x\eta-y\xi, \,\,\,\,\,q_2=x\xi+y\eta$.
If $k_{\op{e}},\, k_{\op{h}},\, k_{\op{f}}$ respectively denote the
number of elliptic, hyperbolic and focus\--focus components, we may
associate the triple $(k_{\op{e}},\, k_{\op{h}},\, k_{\op{f}})$ to a
singularity. The triple is called by Nguy{\^e}n Ti{\^e}n Zung the
\emph{Williamson type of the singularity}.

\begin{theorem}[Eliasson \cite{eliasson,eliasson-these, christophesan}]
  The non\--degenerate critical points of a completely integrable
  system $F \colon M \to \mathbb{R}^n$ are linearizable, i.e. if $m
  \in M$ is a non-degenerate critical point of the completely
  integrable system $F=(f_1,\,\ldots,f_n): M \rightarrow \mathbb{R}^n$
  then there exist local symplectic coordinates $(x_1,\, \ldots,x_n,\,
  \xi_1,\, \ldots,\, \xi_n)$ about $m$, in which $m$ is represented as
  $(0,\,\ldots,\, 0)$, such that $\{f_i,\,q_j\}=0$, for all indices
  $i,\,j$, where we have the following possibilities for the
  components $q_1,\,\ldots,\,q_n$, each of which is defined on a small
  neighborhood of $(0,\,\ldots,\,0)$ in $\mathbb{R}^n$:
  \begin{itemize}
  \item[{\rm (i)}] Elliptic component: $q_j = (x_j^2 + \xi_j^2)/2$,
    where $j$ may take any value $1 \le j \le n$.
  \item[{\rm (ii)}] Hyperbolic component: $q_j = x_j \xi_j$, where $j$
    may take any value $1 \le j \le n$.
  \item[{\rm (iii)}] Focus\--focus component: $q_{j-1}=x_{j-1}\,
    \xi_{j} - x_{j}\, \xi_{j-1}$ and $q_{j} =x_{j-1}\, \xi_{j-1}
    +x_{j}\, \xi_{j}$ where $j$ may take any value $2 \le j \le n-1$
    (note that this component appears as ``pairs'').
  \item[{\rm (iv)}] Non\--singular component: $q_{j} = \xi_{j}$, where
    $j$ may take any value $1 \le j \le n$.
  \end{itemize}
  Moreover if $m$ does not have any hyperbolic block, then the system
  of commuting equations $\{f_i,\,q_j\}=0$, for all indices $i,\,j$,
  may be replaced by the single equation
 $$
 (F-F(m))\circ \varphi = g \circ (q_1,\, \ldots,\,q_n),
 $$ 
 where $\varphi=(x_1,\, \ldots,x_n,\, \xi_1,\, \ldots,\, \xi_n)^{-1}$
 and $g$ is a diffeomorphism from a small neighborhood of the origin
 in $\mathbb{R}^n$ into another such neighborhood, such that $g(0,\,
 \ldots,\,0)=(0,\,\ldots,\,0)$.
\end{theorem}

If the dimension of $M$ is $4$ and $F$ has no hyperbolic singularities
-- which is the case which is most important to us in this paper -- we
have the following possibilities for the map $(q_1,\,q_2)$, depending
on the rank of the critical point:
\begin{itemize}
\item[{\rm (1)}] if $m$ is a critical point of $F$ of rank zero, then
  $q_j$ is one of
  \begin{itemize}
  \item[{\rm (i)}] $q_1 = (x_1^2 + \xi_1^2)/2$ and $q_2 = (x_2^2 +
    \xi_2^2)/2$.
  \item[{\rm (ii)}] $q_1=x_1\xi_2 - x_2\xi_1$ and $q_2 =x_1\xi_1
    +x_2\xi_2$; \,\, on the other hand,
  \end{itemize}
\item[(2)] if $m$ is a critical point of $F$ of rank one, then
  \begin{itemize}
  \item[{\rm (iii)}] $q_1 = (x_1^2 + \xi_1^2)/2$ and $q_2 = \xi_2$.
  \end{itemize}
\end{itemize}
In this case, a non\--degenerate critical point is respectively called
\emph{elliptic\--elliptic, transversally\--elliptic} or
\emph{focus\--focus} if both components $q_1,\, q_2$ are of elliptic
type, one component is of elliptic type and the other component is
$\xi$, or $q_1,\,q_2$ together correspond to a focus\--focus
component.

The analytic case of Eliasson's theorem was proved by
R{\"u}{\ss}mann~\cite{russmann} for two degrees of freedom systems
($2n=4$) and by Vey~\cite{vey} in any dimension. In the
$\op{C}^{\infty}$ category the \emph{lemme de Morse isochore} of Colin
de Verdi{\`e}re and Vey~\cite{colin-vey} implies Eliasson's result for
one degree of freedom systems. Eliasson's proof of the general case
was somewhat loose at a crucial step, but recently this has been
clarified~\cite{san-miranda}.

\subsubsection{Degenerate critical points} \label{deg:sec}

Degenerate critical points appear in many applications, i.e. rigid
body dynamics. The study of non\--degenerate critical points of
integrable systems is difficult, and little is known in general.  A
few particular situations are relatively understood.  For analytic
systems with one degree of freedom a more concrete method is presented
in \cite{10}.  A general linearization result in the analytic category
is given in \cite{68}. Another result may be found in
\cite{zungmore}. Further study of degenerate singularities may be
found in Kalashnikov \cite{Ka1998} where a semiglobal topological
classification of stable degenerate singularities of corank 1 for
systems with two degrees of freedom is given,
Bolsinov\--Fomenko\--Richter \cite{BoFoRi2000} where it was shown how
semiglobal topological invariants of degenerate singularities can be
used to describe global topological invariants of integrable systems
with two degrees of freedom and Nekhoroshev\--Sadovskii\--Zhilinskii
\cite{NeSaZh2006}, where the so called fractional monodromy phenomenon
is explained via topological properties of degenerate singularities
corresponding to higher order resonances. The best context to approach
these singularities is probably algebraic geometry, and we hope that
this article may bring some additional interactions between algebraic
geometers and specialists on integrable systems; interactions have
began to develop in the context of mirror symmetry, where the study of
singular Lagrangian fibrations is relevant \cite{grs1, grs2, grs3,
  grs4, castano-1, castano0, castano}. From the viewpoint of algebraic
singularity theory, several interesting results concerning the
deformation complex of Lagrangian varieties have been discovered
recently~\cite{sevenheck,garay}.

\emph{Throughout this paper, and unless otherwise stated, we assume
  that all singularities are non\--degenerate.}

\section{Semiglobal theory of completely integrable
  systems} \label{sec:semilocalsystems}

If one aims at understanding the geometry of a completely integrable
foliation or its microlocal analysis, the \emph{semiglobal} aspect is
probably the most fundamental. The terminology semiglobal refers to
anything that deals with an invariant neighborhood of a leaf of the
foliation. This semiglobal study is what allows for instance the
construction of quasimodes associated to a Lagrangian
submanifold. Sometimes semiglobal merely reduces to local, when the
leaf under consideration is a critical point with only elliptic
blocks.

\subsection{Regular fibers}

The analysis of neighborhoods of regular fibers, based on the so
called \aangles{} theorem (also known as the action\--angle theorem)
is now routine and fully illustrated in the literature, for classical
aspects as well as for quantum ones. It is the foundation of the whole
modern theory of completely integrable systems in the spirit of
Duistermaat's seminal article~\cite{duistermaat}, but also of KAM-type
perturbation theorems.

The microlocal analysis of action-angle variables starts with the work
of Colin de Verdi{\`e}re~\cite{colinII}, followed in the $\h$
semiclassical theory by Charbonnel~\cite{charbonnel} and more recently
by the second author and various articles by Zelditch, Toth, Popov,
Sj{\"o}strand and many others. The case of compact symplectic
manifolds has recently started, using the theory of Toeplitz
operators~\cite{charles-bs}.

\begin{figure}[htbp]
  \begin{center}
    \includegraphics[width=10cm]{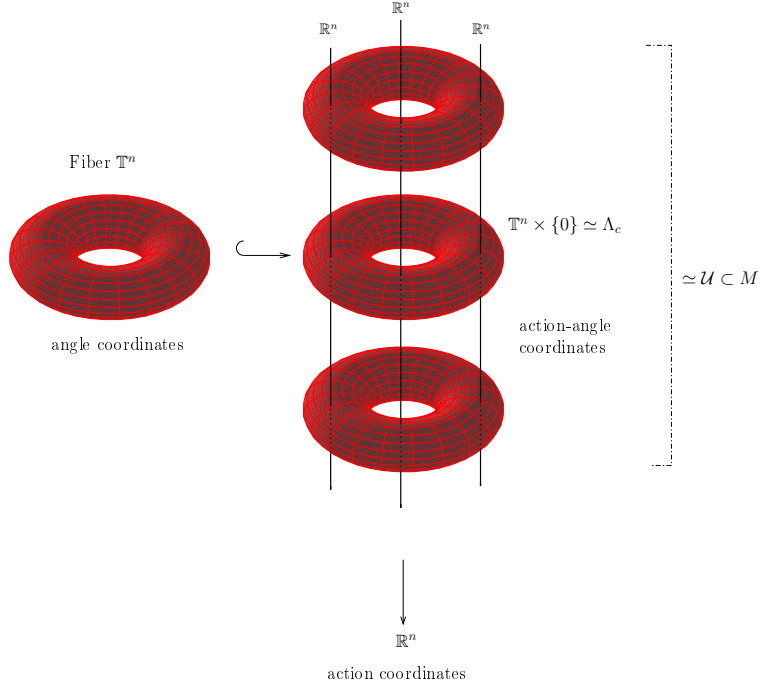}
    \caption{According to the Arnold\--Liouville\--Mineur theorem a
      tubular neighborhood $\mathcal{U}$ of a regular fiber
      $\Lambda_c$ embeds symplectically into
      $\op{T}^*(\mathbb{T}^n)\simeq \mathbb{T}^n \times
      \mathbb{R}^n$.}
    \label{figure1}
  \end{center}
\end{figure}

Let $(f_1,\dots,f_n)$ be an integrable system on a symplectic manifold
$M$. For the remainder of this article we shall assume the momentum
map $F$ to be \emph{proper}, in which case all fibers are compact. Let
$c$ be a regular value of $F$. If we restrict to an adequate invariant
open set, we can always assume that the fibers of $F$ are
connected. Let $\Lambda_c:=F^{-1}(c)$. The fibers being compact and
parallelizable (by means of the vector fields $\mathcal{H}_{f_i}$),
they are tori.  In what follows we identify $\op{T}^*\T^n$ with
$\T^n\times\RM^n$, where $\T=\RM/\ZM \simeq S^1$, equipped with
coordinates $(x_1,\dots,x_n,\xi_1,\dots,\xi_n)$ such that the
canonical Liouville $1$\--form is $\sum_i \xi_i \op{d}\!x_i$.

\begin{theorem}[\aangles\, \cite{mineur,arnold}]
  If $\Lambda_c$ is regular, there exists a local symplectomorphism
  $\chi$ from the cotangent bundle $\op{T}^*\T^n$ of $\T^n$ into $M$
  sending the zero section onto the regular fiber $\Lambda_c$ in such
  a way that that $F \circ \chi= \phy \circ (\xi_1,\dots,\xi_n)$ for
  $\phy$ a local diffeomorphism of $\RM^n$.
\end{theorem}

In this statement, $\chi$ is defined on a neighborhood of the zero
section $\{\xi=0\}$, with values in saturated neighborhood of the torus
$\Lambda_c$; on the other hand $\phy$ is defined in a neighborhood of
the origin in $\RM^n$, and $\phy(0)=F(\Lambda_c)$.

It is important to remark that $\op{d}\!\phy$ is an invariant of the
system since it is determined by periods of periodic trajectories of
the initial system. Regarded as functions on $M$ the $\xi_j$'s are called \emph{action
  variables} of the system for one can find a primitive $\alpha$ of
$\omega$ in a neighborhood of $\Lambda_c$ such that the $\xi_j$'s are
integrals of $\alpha$ on a basis of cycles of $\Lambda_c$ depending
smoothly on $c$.  The coordinates $(x_1,\dots,x_n,\xi_1,\dots,\xi_n)$
are known as \emph{action\--angle variables}. See \cite{miranda} for a
version of the action\--angle theorem in the case of Poisson
manifolds.

\subsection{Singular fibers}
\label{sec:sing-fibres}

This section is devoted to the semi\--global structure of fibers with
non\--degenerate singularities. We are only aware of a very small
number of semi\--global results for degenerate singularities as
mentioned in Section \ref{deg:sec}.

The topological analysis of non-degenerate singular fibers was mainly
initiated by Fomenko~\cite{fomenko} and was successfully expanded by a
number of his students cf.~\cite{bolsinov-fomenko-book}. As far as we
know Lerman and Umanskii \cite{LeUm1988, LeUm1995} were certainly
among the first authors who systematically studied critical points of
Poisson actions on symplectic $4$\--manifolds; their paper
\cite{LeUm1988} is an English translation of their original paper
which was written in Russian and published in a local journal. This
paper is probably the first where focus\--focus singularities are
treated in detail. These works by Lerman and Umanskii have had an
important influence on the Fomenko school. We would also like to mention that
M. Kharlamov appears to be the first author who systematically did a
topological analysis of integrable systems in rigid body dynamics. His
results and methods \cite{Kh1988} were precursors of various aspects
of the mathematical theory of Hamiltonian systems that we discuss in
the present paper (unfortunately \cite{Kh1988} has not been translated
into English yet, but some of the references therein refer to English
papers of the author where his original results may be found).

\begin{figure}[htbp]
  \begin{center}
    \includegraphics[width=10cm]{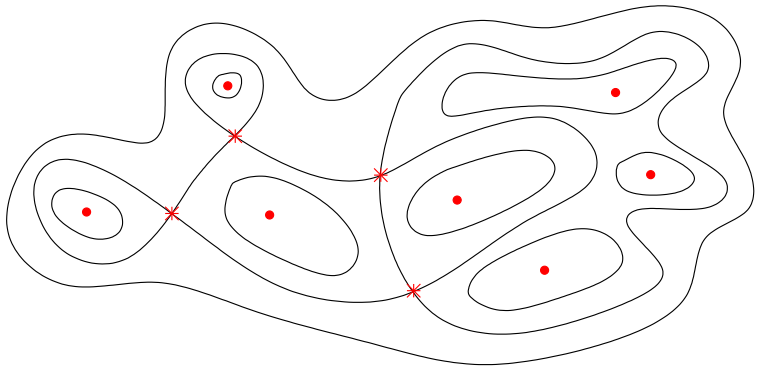}
    \caption{Elliptic and hyperbolic singularities on a surface (the
      figure shows how in the same leaf there may be several
      singularities, though we make the generic assumption that there
      is at most one singularity per fiber). Hyperbolic singularities
      are represented by a red star, elliptic singularities are
      represented by a red dot. Note how nearby a hyperbolic
      singularity the local model looks as in Figure
      \ref{hyperbolic-zoom}. Around an elliptic singularity the leaves
      are concentric circles around the singularity.}
    \label{hyperbolic-global}
  \end{center}
\end{figure}

\subsubsection{Elliptic case} 
Near an elliptic fixed point, the fibers are small tori and are
entirely described by the local normal form, for classical systems as
well as for semiclassical ones (the system is reduced to a set of
uncoupled harmonic oscillators). Therefore we shall not talk about
this type of singularity any further, even if strictly speaking the
semi\--global semiclassical study has not been fully carried out for
\emph{transversally elliptic} singularities. But no particular
difficulties are expected in that case.

\subsubsection{Hyperbolic case}

Just as elliptic blocks, hyperbolic blocks have $1$ degree of freedom
(normal form $q_i=x_i\xi_i$); but they turn out to be more
complicated.  However, in the particular case that $M$ is a surface
there is a classification due to Dufour, Molino and Toulet, which we
present next. In addition to the result, this classification is
interesting to us because it introduces a way to construct symplectic
invariants which is similar to the way symplectic invariants are
constructed for focus\--focus singularities (which will be key to
study semitoric integrable systems in sections \ref{sec:semitoric},
\ref{sec:inv}, \ref{sec:semitoric2}, of this paper). Moreover, using
this classification as a stepping stone, Dufour, Molino and Toulet
gave a global symplectic classification of completely integrable
systems on surfaces, which serves as an introduction to the recent
classification of semitoric integrable systems on symplectic
$4$\--manifolds given later in the paper.

Two one\--degree of freedom completely integrable systems $(M,\,
\omega,\, f_1)$ and $(M,\, \omega',\, f'_1)$ are \emph{isomorphic} if
there exists a symplectomorphism $\chi \colon M \to M'$ and a smooth
map $g$ such that $ \chi^*f'_1=g \circ f_1$.  A (non\--degenerate)
critical point $p$ of $(M,\, \omega,\, f_1)$ is either elliptic or
hyperbollic (there cannot be focus\--focus points).

\begin{figure}[htbp]
  \begin{center}
    \includegraphics[width=12cm]{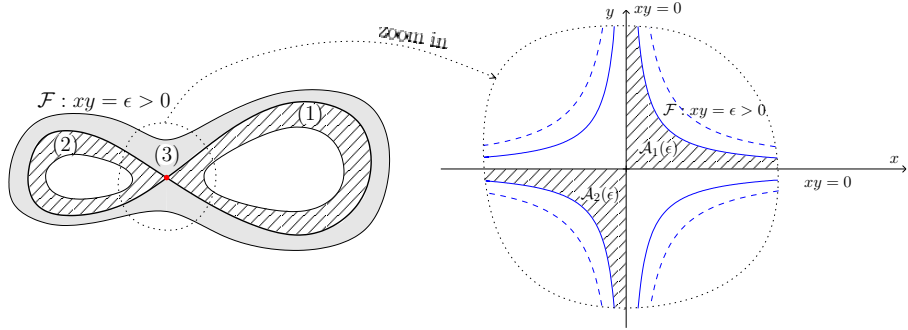}
    \caption{Zoom in around a hyperbolic singularity at the
      intersection of the $x$ and $y$ axes.}
    \label{hyperbolic-zoom}
  \end{center}
\end{figure}

If $p$ is elliptic, there exist local coordinates $(x,\,y)$ and a
function $g$ in a neighborhood $V$ of $p$ such that $f_1=g(x^2+y^2)$
and $\omega=\op{d}\!x \wedge \op{d}\!y$, so geometrically the integral
curves of the Hamiltonian vector field $\mathcal{H}_{f_1}$ generated
by $f_1$ are concentric circles centered at $p$.  If $p$ is
hyperbolic, there exist local coordinates $(x,\,y)$ in a neighborhood
$U$ of $p$ and a function $h$ such that $f_1=h(xy)$ and
$\omega=\op{d}\!x \wedge \op{d}\!y$. In this case the integral curves
of $\mathcal{H}_{f_1}$ are hyperboloid branches
$xy=\textup{constant}$.  One usually calls these integral curves the
\emph{leaves} of the foliation induced by $f_1$.  We make the generic
assumption that our systems may have at most one singularity per leaf
of the induced foliation.

The saturation of $U$ by the foliation has the appearance of an
enlarged figure eight with three components, (1) and (2) corresponding
to $xy>0$, and (3) corresponding to $xy<0$, see Figure
\ref{hyperbolic-zoom} and Figure \ref{hyperbolic-global}, where
therein $\mathcal{F}$ denotes the entire leaf of the foliation
generated by $f_1$ defined in local coordinates near $p$ by
$xy=\epsilon>0$, for some $\epsilon>0$. The area in region (1) between
$\mathcal{F}$ and the figure eight defined locally by $xy=0$ is given
by $ \mathcal{A}_1(\epsilon)=-\epsilon \ln(\epsilon)+h_1(\epsilon) $
for some smooth function $h_1=h_1(\epsilon)$. Similarly the area in
region (2) between $\mathcal{F}$ and the figure eight $xy=0$ is given
by $ \mathcal{A}_2(\epsilon)=-\epsilon\ln(\epsilon)+h_2(\epsilon) $
for some smooth function $h_2=h_2(\epsilon)$, and the area in region
(3) between $\mathcal{F}$ and the figure eight is given by $
\mathcal{A}_3(\epsilon)=2\epsilon \ln |\epsilon|+h(\epsilon) $ for
some smooth function $h=h(\epsilon)$. In addition, for each $q$,
 $$
 h^{(q)}=-h_1^{(q)}(0)+h_2^{(q)}(0).
 $$
 The Taylor series at $0$ of $h_1$ and $h_2$ are symplectic invariants
 of $(M,\, \omega,\, f_1)$, cf. \cite[Proposition 1
 ]{dufour-mol-toul}.
 
 Let $\mathcal{G}$ be the topological quotient of $M$ by the relation
 $a \sim b$ if and only if $a$ and $b$ are in the same leaf of the
 foliation induced by $f_1$. The space $\mathcal{G}$ is called the
 \emph{Reeb graph of $M$}. Let $\pi \colon M \to \mathcal{G}$ be the
 canonical projection map. In this context we call a \emph{regular
   point} the image by $\pi$ of a regular leaf, a \emph{bout} the
 image of an elliptic point, and a \emph{bifurcation point} the image
 of a figure eight leaf. The \emph{edges} are the parts of
 $\mathcal{G}$ contained between two singular points. If $s$ is a
 bifurcation point, then it has three edges, one of which corresponds
 to the leaves extending along the separatrix. We say that this edge
 is the \emph{trunk of $s$}. The two other edges are the
 \emph{branches of $s$}. The graph $\mathcal{G}$ is provided with the
 measure $\mu$, the image by $\pi$ of the measure defined by the
 symplectic form $\omega$ on $M$.

 Let $p \in M$ be a hyperbolic point and let $(x,\,y)$ be the
 aforementioned local coordinates in a neighborhood of $p$.  We define
 a function $\epsilon=xy$ in a neighborhood of the corresponding
 bifurcation point $s$, such that $xy>0$ on the branches of $s$ and
 $xy<0$ on the trunk of $s$.  The expression of $\mu$ in this
 neighborhood is:
 \begin{eqnarray} \label{li} \left\{
     \begin{array}{ccl}
       \op{d}\!\mu_i(\epsilon)&=&(\op{ln}(\epsilon)+g_i(\epsilon))\op{d}\!\epsilon\,\,\, 
       \textup{on each branch}\,\, i=1,\,2
       \\
       \op{d}\!\mu(\epsilon)&=&(2\op{ln}|\epsilon|+g(\epsilon))\op{d}\!\epsilon \,\,\, \textup{on the
         trunk, with} \,\, g,\,g_1,\,g_2 \\
       &&\textup{smooth functions satisfying, for each}\,\, q,  \\
       g^{(q)}(0)&=&((g_1^{(q)}(0)+g_2^{(q)}(0)). 
     \end{array}
   \right.
 \end{eqnarray}

 It follows from \cite[Proposition 1]{dufour-mol-toul} stated above
 that the Taylor series at $0$ of $g_1,\,g_2$ are invariants of
 $(\mathcal{G},\, \mu)$.
 
 \begin{definition}[D{\'e}finition 1 in
   \cite{dufour-mol-toul}] \label{def:1} Let $\mathcal{G}$ be a
   topological $1$\--complex whose vertices have degrees $1$ or
   $3$. For each degree $3$ vertex $s$, which one calls a
   \emph{bifurcation point}, one distinguishes an edge and calls it
   the \emph{trunk of} $s$; the two others are the \emph{branches of
     $s$}. We provide $\mathcal{G}$ with an atlas of the following
   type:
   \begin{itemize}
   \item Outside of the bifurcation points it is a classical atlas of
     a manifold with boundary of dimension $1$.
   \item In a neighborhood of each bifurcation point $s$, there exists
     an open set $V$ and a continuous map $\varphi \colon V \to
     (-\epsilon, \, \epsilon)$, $\epsilon>0$, with $\varphi(s)=0$ and
     such that, if $T$ is the trunk of $s$ and $B_1,\, B_2$ are the
     branches of $s$, $\varphi|_{B_i}$ is bijective on $[0,\,
     \epsilon)$, $i=1,\,2$ and $\varphi|_T$ is bijective on
     $(-\epsilon,\, 0]$. We require that the changes of charts are
     smooth on each part $T\cup B_i$, $i=1,\,2$.
   \end{itemize}

   The topological $1$\--complex $\mathcal{G}$ is provided with a
   measure given by a non\--zero density, smooth on each edge, and
   such that for each vertex $s$ of degree $3$, there exists a chart
   $\varphi$ at $s$ in which the measure $\mu$ is written as in
   \ref{li}.  We denote by $(\mathcal{G},\, \mathcal{D},\, \mu)$ such
   a graph provided with its smooth structure and its measure, and we
   call it an \emph{affine Reeb graph}. An \emph{isomorphism} of such
   a graph is a bijection preserving the corresponding smooth
   structure and measure.
 \end{definition}
 
 One can show \cite[Lemme~2]{dufour-mol-toul} that if the measure
 $\mu$ is written in another chart $\tilde{\varphi}$ in a neighborhood
 of $s$,
 \begin{eqnarray}
   \op{d}\!\mu_i(\epsilon)&=&(\op{ln}(\epsilon)+\tilde{g}_i(\epsilon))\op{d}\!\epsilon \,\,\,\,\,\,
   \textup{on each branch}\,\, B_i \,\, \textup{of}\,\,s \nonumber \\
   \op{d}\!\mu(\epsilon)&=&(2\op{ln}|\epsilon|+\tilde{g}(\epsilon))\op{d}\!\epsilon \,\,\,\,\,\,
   \textup{on the trunk of}\,\,s, \nonumber 
 \end{eqnarray}
 then the functions $g_i$ and $\tilde{g}_i$ have the same Taylor
 series at the origin (hence the Taylor series of $g$ and $\tilde{g}$
 are equal), which shows that the Taylor series of the functions $g_i$
 give invariants for the bifurcation points.

 Let $\mathcal{G}$ be a combinatorial graph with vertices of degree
 $1$ of of degree $3$. For each vertex $s$ of degree $3$, one
 distinguishes in the same fashion as in Definition \ref{def:1} the
 trunk and the branches of $s$, and one associates to each a sequence
 of real numbers. In addition, to each edge one associates a positive
 real number, called its \emph{length}.  Such a graph is called a
 \emph{weighted Reeb graph}. To each affine Reeb graph one naturally
 associated a weighted Reeb graph, the sequence of numbers associated
 to the branches corresponding to coefficients of the Taylor series of
 the functions $g_i$. The lengths of the edges are given by their
 measure.
 
 These considerations tell us the first part of the the following
 beautiful classification theorem.
  
  \begin{theorem}[Dufour\--Molino\--Toulet \cite{dufour-mol-toul}]
    One can associate to a triplet $(M,\, \omega,\, f_1)$ an affine
    Reeb graph $(\mathcal{G}, \, \mathcal{D},\, \mu)$, which is unique
    up to isomorphisms, and to such an affine Reeb graph a weighted
    Reeb graph, unique up to isomorphisms. Conversely, every weighted
    Reeb graph is the graph associated to an affine Reeb graph, unique
    up to isomorphisms, and every affine Reeb graph is the Reeb graph
    associated to a triplet $(M, \, \omega,\, f_1)$, unique up to
    isomorphisms.
  \end{theorem}

  The higher dimensional case will not be treated in general in the
  present paper, as will assume our systems do not have hyperbolic
  singularities.

  \subsubsection{Focus-focus case}

  \emph{Unless otherwise stated, for the remaining of the paper we
    focus on the case when the symplectic manifold $M$ is $4$
    dimensional.}  Eliasson's theorem gives the local structure of
  focus-focus singularities. Several people noticed in the years
  1996-1997 that this was enough to determine the \emph{monodromy} of
  the foliation around the singular fiber.  Actually this local
  structure is a starting point for understanding much more: the
  semiglobal classification of a singular fiber of focus-focus
  type. Unlike monodromy which is a topological invariant, already
  observed in torus fibrations without Hamiltonian structure, the
  semiglobal classification involves purely symplectic invariants.  We
  proceed to describe this semiglobal classification, which
  complements Eliasson's theorem, in two steps.

\begin{figure}[htbp]
    \begin{center}
      \includegraphics[width=6cm]{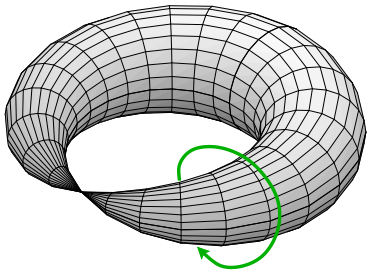}
      \caption{Focus\--focus singularity and vanishing cycle for the
        pinched torus.}
      \label{fig:pinched}
    \end{center}
  \end{figure}

  (a) \emph{Application of Eliasson's theorem.}  Let $F=(f_1,\,f_2)$
  be a completely integrable system with two degrees of freedom on a
  4-dimensional symplectic manifold $M$.  Let $\mathcal{F}$ be the
  \emph{associated singular foliation} to the completely integrable
  system $F=(f_1,\, f_2)$, the leaves of which are by definition the
  connected components of the fibers $F^{-1}(c)$ of $F \colon M \to
  \R^2$.  Let $m$ be a critical point of focus-focus type.  We assume
  for simplicity that $F(m)=0$, and that the (compact, connected)
  fiber $\Lambda_0:=F^{-1}(0)$ does not contain other critical points.
  One can show that $\Lambda_0$ is a ``pinched''
  torus\footnote{Lagrangian immersion of a sphere $S^2$ with a
    transversal double point} surrounded by regular fibers which are
  standard $2$\--tori, see Figure \ref{fig:pinched}.  What are the
  semi\--global invariants associated to this singular fibration?

  One of the major characteristics of focus\--focus singularities is
  the existence of a \emph{Hamiltonian action of $S^1$} that commutes
  with the flow of the system, in a neighborhood of the singular fiber
  that contains $m$.  By Eliasson's theorem \cite{eliasson-these}
  there exist symplectic coordinates $(x,\, y,\, \xi,\,\eta)$ in a
  neighborhood $U$ around $m$ in which $(q_1,\,q_2)$, given by
  \begin{equation}
    q_1=x\eta-y\xi, \,\,   q_2=x\xi+y\eta
    \label{equ:cartan}
  \end{equation}
  is a momentum map for the foliation $\mathcal{F}$; here the critical
  point $m$ corresponds to coordinates $(0,\,0,\,0,\,0)$.  Fix $A'\in
  \Lambda_0\cap (U\setminus\{m\})$ and let $\Sigma$ denote a small
  2\--dimensional surface transversal to $\mathcal{F}$ at the point
  $A'$.
  
  Since the Liouville foliation in a small neighborhood of $\Sigma$ is
  regular for both $F$ and $q=(q_1,\,q_2)$, there is a diffeomorphism
  $\varphi$ from a neighborhood $U'$ of $F(A')\in\R^2$ into a
  neighborhood of the origin in $\R^2$ such that $q=\varphi \circ
  {F}$. Thus there exists a smooth momentum map $\Phi=\varphi
  \circ{F}$ for the foliation, defined on a neighborhood
  $\Omega=F^{-1}(U')$ of $\Lambda_0$, which agrees with $q$ on $U$.

  Write $\Phi:=(H_1,\,H_2)$ and $\Lambda_z:=\Phi^{-1}(c)$.  Note that
  $ \Lambda_0=\mathcal{F}_m.  $ It follows from~(\ref{equ:cartan})
  that near $m$ the $H_1$\--orbits must be periodic of primitive
  period $2\pi$, whereas the vector field $\mathcal{H}_{H_2}$ is
  hyperbolic with a local stable manifold (the $(\xi,\eta)$-plane)
  transversal to its local unstable manifold (the
  $(x,y)$-plane). Moreover, $\mathcal{H}_{H_2}$ is \emph{radial},
  meaning that the flows tending towards the origin do not spiral on
  the local (un)stable manifolds.
  
  (b) \emph{Symplectic semi\--global classification of focus\--focus
    point.}  Suppose that $A \in\Lambda_c$ for some regular value $c$.
  Let $\tau_2(c)>0$ be the time it takes the Hamiltonian flow
  associated with $H_2$ leaving from $A$ to meet the Hamiltonian flow
  associated with $H_1$ which passes through $A$. The existence of
  $\tau_2$ is ensured by the fact that the flow of $H_2$ is a
  quasiperiodic motion always transversal to the $S^1$-orbits
  generated by $H_1$.

  Let $\tau_1(c)\in\R/2\pi\Z$ the time that it takes to go from this
  intersection point back to $A$, closing the trajectory.
  
  The commutativity of the flows ensure that $\tau_1(c)$ and
  $\tau_2(c)$ do not depend on the initial point $A$. Indeed, if
  $\phy_1^t$, $\phy_2^t$ denote the hamiltonian flows of $H_1$ and
  $H_2$, respectively, then we have
  $A=\phy_1^{\tau_1(c)}\circ\phy_2^{\tau_2(c)}(A)$. If
  $\tilde{A}\in\Lambda_c$ is close to $A$, the fact that $\op{d}\!H_1$
  and $\op{d}\!H_2$ are independent near $A$ implies that the
  corresponding hamiltonian action is locally free on $\Lambda_c$~:
  there exist small times $(t_1,t_2)\in\RM^2$ such that
\[
\tilde{A} = \phy_1^{t_1}\circ \phy_2^{t_2}(A).
\]
Thus we may write
\[
\tilde{A} = \phy_1^{t_1+\tau_1(c)}\circ\phy_2^{t_2+\tau_2(c)}(A) =
\phy_1^{\tau_1(c)}\circ\phy_2^{\tau_2(c)} (\tilde{A}).
\]
This shows that the times $\tau_1$, $\tau_2$ that would be obtained
starting from $\tilde{A}$ are the same as those we obtained starting
from $A$.

  Write $c=(c_1,\,c_2)=c_1+\op{i}c_2$ ($c_1,\, c_2 \in \R$), and let
  $\op{ln} c$ be a fixed determination of the logarithmic function on
  the complex plane. Let
  \begin{equation}
    \left\{
      \begin{array}{ccl}
        \sigma_1(c) & = & \tau_1(c)-\Im(\op{ln} c) \\
        \sigma_2(c) & = & \tau_2(c)+\Re(\op{ln} c),
      \end{array}
    \right.
    \label{equ:sigma}
  \end{equation}
  where $\Re$ and $\Im$ respectively stand for the real and imaginary
  parts of a complex number.  V\~u Ng\d oc proved in
  \cite[Proposition~\,3.1]{vungoc0} that $\sigma_1$ and $\sigma_2$ extend to
  smooth and single\--valued functions in a neighborhood of $0$ and
  that the differential 1\--form $\sigma:=\sigma_1\,
  \DD{}c_1+\sigma_2\, \DD{}c_2 $ is closed.  Notice that if follows
  from the smoothness of $\sigma_1$ that one may choose the lift of
  $\tau_1$ to $\R$ such that $\sigma_1(0)\in[0,\,2\pi)$. This is the
  convention used throughout.

  Following \cite[Def.~3.1]{vungoc0} , let $S$ be the unique smooth
  function defined around $0\in\R^2$ such that
  \begin{eqnarray}
    \DD{}S=\sigma,\,\, \,\, S(0)=0.
   \end{eqnarray}
  The Taylor expansion of $S$ at $(0,\,0)$ is denoted by $(S)^\infty$.

  Loosely speaking, one of the components of the system is indeed
  $2\pi$\--periodic, but the other one generates an arbitrary flow
  which turns indefinitely around the focus\--focus singularity,
  deviating from periodic behavior in a logarithmic fashion, up to a
  certain error term; this deviation from being logarithmic is the
  symplectic invariant $(S)^{\infty}$.

  \begin{theorem}[V\~ u Ng\d oc \cite{vungoc0}]
    \label{theo:invariants}
    The Taylor series expansion $(S)^\infty$ is well-defined (it does
    not depend on the choice of Eliasson's local chart) and it
    classifies the singular foliation in a neighborhood of $\Lambda_0$
    in the sense that another system has the same Taylor series
    invariant near a focus\--focus singularity if and only if there is
    a symplectomorphism which takes foliated a neighborhood of the
    singular fiber to a foliated neighborhood of the singular fiber
    preserving the leaves of the foliation and sending the singular
    fiber to the singular fiber.

    Moreover, if $S$ is any formal series in $\RM\formel{X,Y}$ with
    $X$\--coefficient in $[0,\, 2 \pi)$ and without constant term,
    then there exists a singular foliation of focus-focus type whose
    Taylor series expansion is $S$.

  \end{theorem}

  The fact that two focus-focus fibrations are always semiglobally
  \emph{topologically} conjugate was already proved by
  Lerman\--Umanskii and Nguy{\^e}n Ti{\^e}n Zung~\cite{LeUm1988,
    zung-I}, who introduced various topological notions of
  equivalence.

  $S$ can be interpreted as a \emph{regularized} (or
  \emph{desingularized}) \emph{action}. Indeed if $\gamma_z$ is the
  loop on $\Lambda_z$ defined just as in the description of $\tau_j$
  above, and if $\alpha$ is a semiglobal primitive of the symplectic
  form $\omega$, let $\mathcal{A}(z)=\int_{\gamma_c}\alpha$; then
  $$
  S(z) = \mathcal{A}(z) - \mathcal{A}(0) + \re(z\ln z - z).
  $$
  
  \subsection{Example}
  
  One can check that the singularities of the coupled spin--oscillator
  $S^2 \times \R^2$ model mentioned in the Section \ref{sec:intro}
  (equipped with the product symplectic form or the standard area
  forms) are non\--degenerate and of elliptic\--elliptic,
  transversally\--elliptic or focus\--focus type.
    
  It has exactly one focus\--focus singularity at the ``North Pole''
  $((0,\,0,\,1),\,(0,\,0)) \in S^2 \times \R^2$ and one
  elliptic\--elliptic singularity at the ``South Pole''
  $((0,\,0,\,-1),\,(0,\,0))$.  Let us parametriza the singular fiber
  $\Lambda_0:=F^{-1}(1,\,0)$. This singular fiber $\Lambda_0$
  corresponds to the system of equations $J=1$ and $H=0$, which
  explicitly is given by system of two nonlinear equations
  $J=(u^2+v^2)/2 + z=0$ and $H= \frac{1}{2} (ux+vy)=0$ on the
  coordinates $(x,\,y,\,z,\, u,\,v)$.

  In order to solve this system of equations one introduces polar
  coordinates $u+\op{i}v=r\op{e}^{\op{i}t}$ and $
  x+\op{i}y=\rho\op{e}^{\op{i}\theta}$ where recall that the
  $2$\--sphere $S^2 \subset \mathbb{R}^3$ is equipped with coordinates
  $(x,\,y,\,z)$, and $\R^2$ is equipped with coordinates $(u,\, v)$.
  For $\epsilon=\pm 1$, we consider the mapping
  $$
  S_\epsilon : [-1,\, 1]\times \R/2\pi\Z \to \R^2\times S^2
  $$
  given by the formula
$$
S_\epsilon(p) = (r(p)\,\op{e}^{\op{i}t(p)},
\,(\rho(p)\,\op{e}^{\op{i}\theta(p),\, z(p)}))
 $$
 where $p=(\tilde{z}, \, \tilde{\theta})\in [-1,1]\times[0,2\pi)$ and
 \begin{equation} \nonumber
   \begin{cases}
     r(p) = \sqrt{2(1-\tilde{z})}\\
     t(p) = \tilde{\theta} + \epsilon\frac{\pi}{2}\\
     \rho(p)= \sqrt{1-\tilde{z}^2}\\
     \theta(p) = \tilde{\theta}\\
     z(p) = \tilde{z}.
   \end{cases}
 \end{equation}
  
 Then the map $S_{\epsilon}$, where $\epsilon=\pm1$, is continuous and
 $S_{\epsilon}$ restricted to $(-1,\,1) \times \R/2\pi\Z $ is a
 diffeomorphism onto its image.  If we let
 $\Lambda_0^{\epsilon}:=S_{\epsilon}([-1,\, 1]\times \R/2\pi\Z )$,
 then $\Lambda_0^1\cup \Lambda_0^2=\Lambda_0$ and
$$\Lambda_0^1 \cap \Lambda_0^2=\Big( \{(0,\,0)\} \times \{(1,\,0,\,0)\} \Big)\cup 
\Big(C_2 \times \{(0,\,0,\,-1)\} \Big),$$ where $C_2$ denotes the
circle of radius $2$ centered at $(0,\,0)$ in $\R^2$.  Moreover,
$S_{\epsilon}$ restricted to $(-1,\,1) \times \R/2\pi\Z $ is a smooth
Lagrangian embedding into $\R^2 \times S^2$.  The singular fiber
$\Lambda_0$ consists of two sheets glued along a point and a circle;
topologically $\Lambda_0$ is a pinched torus, i.e.  a
$2$\--dimensional torus $S^1 \times S^1$ in which one circle $\{p\}
\times S^1$ is contracted to a point (which is of course not a a
smooth manifold at the point which comes from the contracting
circle). The statements correspond to \cite[Proposition 2.8]{example}.

It was proven in \cite[Theorem 1.1]{example} that the linear deviation
from exhibiting logarithmic in a saturated neighborhood of the
focus\--focus singularity is given by the linear map $L \colon \R^2
\to \R$ with expression
$$
L(X,\, Y)=5\,\,\op{ln}2 \, X+ \frac{\pi}{2}\, Y.
$$
In other words, we have an equality
$$
(S(X,\,Y))^{\infty}=L(X,\,Y) +\mathcal{O}((X,\,Y)^2).
$$
This computation is involved and uses some deep formulas from
microlocal analysis proven in the late 1990s. At the time of writing
this paper we do not have an strategy to compute the higher order
terms of the Taylor series invariant.
 
\subsection{Applications}
Theorem~\ref{theo:invariants} leads to a number of applications, which
although are outside the scope of this paper, we briefly note.  One
can for instance exploit the fact that the set of symplectic
equivalence classes of these foliations acquires a vector space
structure. That is what Symington does in~\cite{symington-gen} to show
that neighborhoods of focus-focus fibers are always symplectomorphic
(after forgetting the foliation, of course). For this one introduces
functions $S_0$ and $S_1$ whose Taylor expansions give the invariants
of the two foliations, and constructs a ``path of foliations'' by
interpolating between $S_0$ and $S_1$. Then a Moser type argument
yields the result (since the symplectic forms are cohomologous).

The theorem is also useful for doing calculations in a neighborhood of
the fiber. For instance it is possible in this way to determine the
validity of non-degeneracy conditions that appear in KAM type
theorems\footnote{A nice discussion of these various conditions can be
  found in~\cite{roy-conditions}}, for a perturbation of a completely
integrable system with a focus-focus singularity (see also
\cite{zung-kolmogorov}).

\begin{theorem}[Dullin\--V\~ u Ng\d oc \cite{san-dullin}]
  Let $H$ be a completely integrable Hamiltonian with a loxodromic
  singularity at the origin (i.e. $H$ admits a singular Lagrangian
  foliation of focus-focus type at the origin). Then Kolmogorov's
  non\--degeneracy condition is fulfilled on all tori close to the
  critical fiber, and the ``isoenergetic turning frequencies''
  condition is fulfilled except on a $1$\--parameter family of tori
  corresponding to a curve through the origin in the image of the
  momentum map which is transversal to the lines of constant energy
  $H$.
\end{theorem}

\subsection{A topological classification}

The present paper is devoted to the symplectic theory of Hamiltonian
integrable systems. A large number authors have made contributions to
the topological theory of Hamiltonian integrable systems, in
particular Fomenko and his students, and Nguy{\^e}n Ti{\^e}n Zung. In
this section we briefly present a classification result due to Zung,
which holds in any dimension, and for all so called topologically
stable (non\--degenerate) singularities.  For precise statements we
refer to Zung \cite[Section 7]{zung-I}.

Let $\mathcal{F}$ be a singular leaf (fiber corresponding to a
non\--degenerate singularity) of an integrable system. In what
follows, a \emph{tubular neighborhood} $\mathcal{U}(\mathcal{F})$ of
$\mathcal{F}$ means an appropriately chosen sufficiently small
saturated tubular neighborhood. We denote by
$(\mathcal{U}(\mathcal{F}), \, \mathcal{L})$ the Lagrangian foliation
in a tubular neighborhood $\mathcal{U}(\mathcal{F})$ of
$\mathcal{F}$. The leaf $\mathcal{F}$ is a deformation retract of
$\mathcal{U}(\mathcal{F})$.

Let $\mathcal{F}_1,\, \mathcal{F}_2$ be (non\--degenerate) singular
leaves of two integrable systems, of coranks $k_1,\, k_2$,
respectively, and let $(\mathcal{U}(\mathcal{F}_1),\, \mathcal{L}_1)$
and $(\mathcal{U}(\mathcal{F}_2),\, \mathcal{L}_2)$ be the
corresponding Lagrangian foliations. Here a singular leaf is said to
be \emph{ of corank $k$} if $k$ is the maximal corank of the
differential of the momentum map at critical points of the singular
leaf\footnote{In \cite{zung-I}, Zung used the terminology
  ``codimension'' for the corank $k$.}. For instance a focus-focus
singularity in a symplectic 4-manifold has corank $2$. The
\emph{direct product} of these singularities is the singular leaf
$\mathcal{F}= \mathcal{F}_1 \times \mathcal{F}_2$ of corank $k_1+k_2$,
with the associated Lagrangian foliation
$$(\mathcal{U}(\mathcal{F}),\, \mathcal{L}):=(\mathcal{U}(\mathcal{F}_1),\, \mathcal{L}_1)
\times (\mathcal{U}(\mathcal{F}_2),\, \mathcal{L}_2).$$

A (non\--degenerate) singularity (or singular pair)
$(\mathcal{U}(\mathcal{F}),\, \mathcal{L})$ of corank $k$ and
Williamson type $(k_{\op{e}},\, k_{\op{h}},\, k_{\op{f}})$ of an
integrable system with $n$ degrees of freedom is called \emph{of
  direct product type topologically} if it is homeomorphic, together
with the Lagrangian foliation, to the direct product
\begin{eqnarray}
  (\mathcal{U}(\mathbb{T}^{n-k}),\, \mathcal{L}_{\op{r}})
  \times (\op{P}^2(\mathcal{F}_1),\, \mathcal{L}_1) \times \ldots \times
  (\op{P}^2(\mathcal{F}_{k_{\op{e}}+k_{\op{h}}}),\, \mathcal{L}_{k_{\op{e}}+k_{\op{h}}}) \times \nonumber \\
  (\op{P}^4(\mathcal{F}'_1),\, \mathcal{L}'_1) \times \ldots \times
  (\op{P}^4(\mathcal{F}'_{k_{\op{f}}}),\, \mathcal{L}'_{k_{\op{f}}}) \label{for:zung}
\end{eqnarray}
where:
\begin{itemize}
\item the tuple $(\mathcal{U}(\mathbb{T}^{n-k},\,
  \mathcal{L}_{\op{r}})$ denotes the Lagrangian foliation in a tubular
  neighborhood of a regular $(n-k)$\--dimensional torus of an
  integrable system with $n-k$ degrees of freedom;
\item the tuple $(\op{P}^2(\mathcal{F}_i),\, \mathcal{L}_i)$, $1 \le i
  \le k_{\op{e}}+k_{\op{h}}$, denotes a codimension $1$
  (non\--degenerate) surface singularity (i.e. a singularity of an
  integrable system with one degree of freedom);
\item the tuple $(\op{P}^4(\mathcal{F}_i),\, \mathcal{L}'_i)$, $1 \le
  i \le k_{\op{f}}$, denotes a focus\--focus singularity of an
  integrable system with two degrees of freedom.
\end{itemize}

In this case we have
$k=k_{\textup{e}}+k_{\textup{h}}+2k_{\textup{f}}$. A
(non\--degenerate) singularity of an integrable system is \emph{of
  almost direct product type topologically} if a finite covering of it
is homeomorphic, together with the Lagrangian foliation, to a direct
product singularity. Nguy{\^e}n Ti{\^e}n Zung proved \cite[Theorem 7.3]{zung-I} the
following classification result.

\begin{theorem}[Nguy{\^e}n Ti{\^e}n Zung
  \cite{zung-I}] \label{thm:zung} If $(\mathcal{U}(\mathcal{F}),\,
  \mathcal{L})$ is a (non\--degenerate) topologically stable
  singularity of Williamson type $(k_{\op{e}},\, k_{\op{h}},\,
  k_{\op{f}})$ and corank $k$ of an integrable system with $n$ degrees
  of freedom, then it can be written homeomorphically in the form of a
  quotient of a direct product singularity as in (\ref{for:zung}) by
  the free action of a finite group $\Gamma$ which acts
  component\--wise on the product, and acts trivially on the elliptic
  components.
\end{theorem}

As Nguy{\^e}n Ti{\^e}n Zung points out, the decomposition in Theorem
\ref{thm:zung} is in general \emph{not} symplectic.

\section{Introduction to semitoric completely integrable
  systems} \label{sec:semitoric}

For the remainder of this paper we are now going to focus exclusively
on semitoric completely integrable systems with $2$\--degrees of
freedom on $4$\--manifolds; for brevity we will call these simply
``semitoric systems''.  Essentially this means that the system is
half\--toric, and half completely general --- but non-toric
singularities must be isolated --- see Definition \ref{semitoric:def}
for a precise definition. Naturally, hamiltonian toric manifolds form
a strict subclass of semitoric systems.

Semitoric systems form an important class of integrable systems,
commonly found in simple physical models. Indeed, a semitoric system
can be viewed as a Hamiltonian system in the presence of an
$S^1$\--symmetry~\cite{sadovski-zhilinski}.  In our personal opinion,
it is much simpler to understand the integrable system on its whole
rather than writing a theory of Hamiltonian systems on Hamiltonian
$S^1$\--manifolds.

\subsection{Meaning of the integrability condition}

Let us recall what the general definition of an integrable system in
Section \ref{sec:intro} means in dimension $4$. 
In this case an integrable system on $M$ is a pair of real\--valued
smooth functions $J$ and $H$ on $M$, for which the Poisson bracket
$\{J,\,H\}:=\omega(\mathcal{H}_J,\, \mathcal{H}_H)$ identically
vanishes on $M$, and the differentials $\op{d}\!J$, $\op{d}\!H$ are
almost\--everywhere linearly independendent.

Of course, here $(J,\,H) \colon M \to \R^2$ is the analogue of the
momentum map in the case of a torus action.  In some local Darboux
coordinates of $M$, $(x,\, y,\, \xi,\, \eta)$, the symplectic form
$\omega$ is given by $\op{d}\!  \xi \wedge \op{d}\!x +\op{d}\!
\eta\wedge \op{d}\!y$, and the vanishing of the Poisson brackets
$\{J,\,H\}$ amounts to the partial differential equation
  $$
  \frac{\partial J}{\partial \xi} \, \frac{\partial H}{\partial x} -
  \frac{\partial J}{\partial x} \, \frac{\partial H}{\partial \xi} +
  \frac{\partial J}{\partial \eta} \, \frac{\partial H}{\partial y} -
  \frac{\partial J}{\partial y} \, \frac{\partial H}{\partial \eta}
  =0.
  $$
  This condition is equivalent to $J$ being constant along the
  integral curves of $\mathcal{H}_{H}$ (or $H$ being constant along
  the integral curves of $\mathcal{H}_{J}$).

  \subsection{Singularities}

  We introduce the main object of the remaining part of this paper,
  semitoric systems, and explain what singularities can occur in these
  systems.

  \begin{definition} \label{semitoric:def} A {\em semitoric integrable
      system on $M$} is an integrable system for which the component
    $J$ is a proper momentum map for a Hamiltonian circle action on
    $M$, and the associated map $F:=(J,\,H):M\to\R^2$ has only
    non\--degenerate singularities in the sense of Williamson, without
    real\--hyperbolic blocks.
  \end{definition}

  \begin{remark}
   There are examples which come endowed with a Hamiltonian
    $S^1$\--action but do not fit Definition \ref{semitoric:def}. 
    The direct product $S^2 \times S^2$ equipped with the Hamiltonians
    $J:=z_1$ and $H:=x_2y_2$, where $(x_1,\,y_1,\,z_1)$,
    $(x_2,\,y_2,\,z_2)$ are respectively the coordinates on the first
    and second copies of $S^2$, is a non-degenerate system with a
    proper $S^1$-momentum map $J$, but it contains hyperbolic
    singularities.

    Similarly, the spherical pendulum on $\op{T}^*\!S^2$ (mentioned in
    Section \ref{se}) equipped with
    $J(\theta,\varphi,\xi_\theta,\xi_\phy)= \xi_\phy$ and
    \[
    H(\underbrace{\theta,\phy}_{\tiny \textup{sphere}},
    \underbrace{\xi_\theta,\xi_\phy}_{\textup{fiber}})=
    \underbrace{\frac{1}{2}\left(\xi_\theta^2 +
        \frac{1}{\sin^2\theta}\xi_\phy^2\right)}_{\textup{kinetic
        energy}}+ \underbrace{\cos\theta}_{\textup{potential}}
    \]
    is a non-degenerate system that does not have hyperbolic
    singularities. However $J$ is not a proper map.  We are developing
    a theory that deals with this situation \cite{PeRaVu2010}.
  \end{remark}

  Let us spell this definition concretely. The hamiltonian flow
  defined by $J$ is periodic, with fixed period $2\pi$. The flow of
  $H$ is in general not periodic, but of course it is quasiperiodic on
  regular Liouville tori. The properness of $J$ means that the
  preimage by $J$ of a compact set is compact in $M$ (which is
  immediate if $M$ is compact). The non\--degeneracy hypothesis for
  $F$ means that, if $m$ is a critical point of $F$, then there exists
  an invertible 2 by 2 matrix $B$ such that, if we write
  $\tilde{F}=B\circ F,$ one of the situations described in the
  following table holds in some local symplectic coordinates
  $(x,\,y,\, \xi,\, \eta)$ near $m$ in which $m=(0,0,0,0)$ and
  $\omega=\op{d}\!  \xi \wedge \op{d}\!x +\op{d}\! \eta\wedge
  \op{d}\!y$.

  \begin{center}
    \begin{tabular}{|l|l|}
      \hline
      {\tiny {\bf TYPE}} & ${F}:=(H,J)\colon M \to \mathbb{R}^2$  
      {\tiny {\bf IN COORDINATES}} $(x,\,y,\, \xi,\, \eta)$  \\ \hline

      {\tiny Transversally elliptic}  &  { ${ F}=
        (\eta+\mathcal{O}(\eta^2),\, \frac12(x^2+\xi^2)+\mathcal{O}((x,\,
        \xi)^3)$
      } \\ \hline

      {\tiny Elliptic-elliptic} &  { ${F}=
        \frac{1}{2}(x^2+\xi^2,\,y^2+\eta^2)+\mathcal{O}((x,\, \xi,\, y,\,
        \eta)^3)$}  \\ \hline

      {\tiny Focus-focus}   &  { ${F}=(x\eta-y\xi,\, x\xi+y\eta) +\mathcal{O}((x,\, \xi,\, y,
        \,\eta)^3)$}   \\ \hline     
    \end{tabular}
  \end{center}

  In the case of semitoric systems the Williamson of the singularities
  are of the form $(k_{\op{e}},\, 0,\, \op{k}_{\op{f}})$,
  i.e. $k_{\op{h}}=0$.

  Note that this is not a result in a neighborhood of a fiber.
 
  Again, perhaps the simplest non\--compact semitoric integrable
  systems is the coupled spin\--oscillator $S^2\times\R^2$. The
  component $J$ is the momentum map for the Hamiltonian circle action
  on $M$ which rotates simultaneously about the vertical axis of $S^2$
  and about the origin of $\R^2$.  By the action\--angle theorem of
  Arnold\--Liouville\--Mineur the regular fibers are $2$\--tori. We
  saw that this system has a unique focus\--focus singularity at
  $(0,0,1,0,0)$ with fiber a pinched torus. The other singular fibers
  are either circles or points.  The authors have studied the
  symplectic and spectral theory of this system in \cite{example}.

  \subsection{Convexity properties: the polygon invariant}
  \label{semitoric:sec}

  This section analyzes to what extent the convexity theorem of Atiyah
  and Guillemin\--Sternberg holds in the context of semitoric
  completely integrable systems. The second author proved in
  \cite{vungoc} that one can meaningfully associate a convex polygonal
  region such a system.

\subsubsection{Bifurcation diagrams}
It is well established in the integrable systems community that the
most simple and natural object, which tells much about the structure
of the integrable system under study, is the so\--called bifurcation
diagram. As a matter of fact, bifurcations diagrams may be defined in
great generality as follows.  Let $M$ and $N$ be smooth
manifolds. Recall that a smooth map $f: M\rightarrow N$ is locally
trivial at $n_0 \in f(M)$ if there is an open neighborhood $U \subset
N$ of $n_0$ such that $f ^{-1}(n)$ is a smooth submanifold of $M$ for
each $n\in U$ and there is a smooth map $h: f^{-1}(U)\rightarrow f
^{-1}(n_0)$ such that $f\times h:f^{-1}(U)\rightarrow U \times f
^{-1}(n_0)$ is a diffeomorphism. The \emph{bifurcation set or
  bifurcation diagram} $\Sigma_f$ consists of all the points of $N$
where $f$ is not locally trivial. Note, in particular, that $h|_{f
  ^{-1}(n)}:f ^{-1}(n) \rightarrow f ^{-1}(n_0)$ is a diffeomorphism
for every $n \in U $. Also, the set of points where $f $ is locally
trivial is open in $N$, so $\Sigma_f $ is a closed subset of $N $.  It
is well known that the set of critical values of $f $ is included in
the bifurcation set (see \cite[Proposition 4.5.1]{AbMa1978}).  In
general, the bifurcation set strictly includes the set of critical
values. This is the case for the momentum-energy map for the two-body
problem \cite[\S9.8]{AbMa1978}. It is well known \cite[Page
340]{AbMa1978} that if $f \colon M \to N$ is a smooth proper map, the
bifurcation set of $f$ is equal to the set of critical values of $f$.

It follows that when the map $F=(J,\,H) \colon M \to \mathbb{R}^2$
that defines the integrable system is a proper map the bifurcation
diagram is equal to the set of critical values of $F$ inside of the
image $F(M)$ of $F$. This is the case for semitoric integrable
systems, since the properness of the component $J$ implies the
properness of $J$. As it turns out, the arrangement of such critical
values is indeed important, but other crucial invariants that are more
subtle and cannot be detected from the bifurcation diagram itself are
needed to understand a semitoric system $F$; we deal with these ones
in Section \ref{sec:inv}.  The authors proved \cite{pelayovungoc1,
  pelayovungoc2} that these invariants are enough to completely
determine a semitoric system up to isomorphisms.

The proof relies on a number of remarkable results by other authors on
integrable systems, including Arnold \cite{arnold}, Atiyah
\cite{atiyah}, Dufour\--Molino \cite{dufour-molino}, Eliasson
\cite{eliasson}, Duistermaat \cite{duistermaat}, Guillemin\--Sternberg
\cite{gs}, Miranda\--Zung \cite{miranda-zung} and V\~u Ng\d oc
\cite{vungoc0, vungoc}. In this section we explain the so called
polygon invariant, which was originally introduced by the second
author in \cite{vungoc}, and can be considered an analogue (for
completely integrable semitoric systems) of the convex polytope that
appears in the Atiyah\--Guillemin\--Sternberg convexity theorem (in
the context of symplectic torus actions on compact manifolds).

\subsubsection{Affine structures}

The plane $\RM^2$ is equipped with its standard affine structure with
origin at $(0,0)$, and its standard orientation.  Let
$\textup{Aff}(2,\RM):=\textup{GL}(2,\RM)\ltimes\RM^2$ be the group of
affine transformations of $\RM^2$. Let
$\textup{Aff}(2,\ZM):=\textup{GL}(2,\ZM)\ltimes\RM^2$ be the subgroup
of \emph{integral-affine} transformations.  It was proven in
\cite{vungoc} that a semitoric system $(M,\, \omega,\, F:=(J,\,H))$
has finitely many focus\--focus critical values
$c_1,\,\dots,\,c_{m_f}$, that if we write $B:=F(M)$ then the set of
regular values of $F$ is $\op{Int}(B)\setminus\{c_1,\dots,c_{m_f}\}$,
that the boundary of $B$ consists of all images of elliptic
singularities, and that the fibers of $F$ are connected. The integer
$m_f$ was the first invariant that we associated with such a system.
  
\begin{figure}[htbp]
  \begin{center}
    \includegraphics[width=0.3\linewidth]{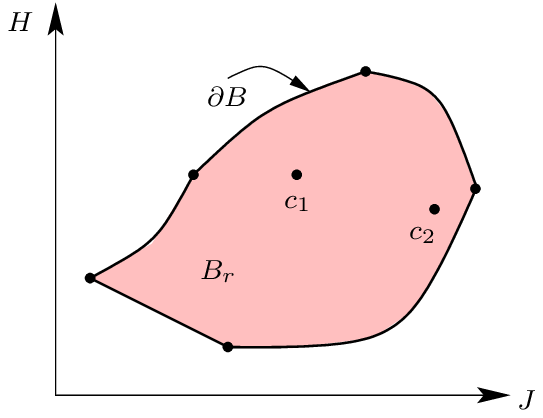}
    \caption{Image under $(J,\,H)$ of $M$.}
    \label{fig:evanescent}
  \end{center}
\end{figure}

Let $\mathfrak{I}$ be the subgroup of $\op{Aff}(2,\,\Z)$ of those
transformations which leave a vertical line invariant, or
equivalently, an element of $\mathfrak{I}$ is a vertical translation
composed with a matrix $T^k$, where $k \in \Z$ and
\begin{eqnarray}
  T^k:=\left(
    \begin{array}{cc}
      1 & 0\\ k & 1
    \end{array}
  \right) \in \op{GL}(2,\, \Z).
  \label{equ:Tk}
\end{eqnarray}

Let $\ell\subset\R^2$ be a vertical line in the plane, not necessarily
through the origin, which splits it into two half\--spaces, and let
$n\in\Z$. Fix an origin in $\ell$.  Let $t^n_{\ell} \colon \R^2 \to
\R^2$ be the identity on the left half\--space, and $T^n$ on the right
half\--space. By definition $t^n_{\ell}$ is piecewise affine.  Let
$\ell_i$ be a vertical line through the focus\--focus value
$c_i=(x_i,\,y_i)$, where $1 \le i \le m_f$, and for any tuple $\vec
n:=(n_1,\,\dots,\,n_{m_f})\in\Z^{m_f}$ we set
$t_{\vec n}:=t^{n_1}_{\ell_1}\circ\, \cdots\, \circ
t^{n_{m_f}}_{\ell_{m_f}}$.
The map $t_{\vec n}$ is piecewise affine.

A \emph{convex polygonal set} $\Delta$ is the intersection in $\RM^2$
of (finitely or infinitely many) closed half\--planes such that on
each compact subset of the intersection there is at most a finite
number of corner points. We say that $\Delta$ is \emph{rational} if
each edge is directed along a vector with rational coefficients.  For
brevity, in this paper we usually write \emph{``polygon''} (or
``convex polygon'') instead of \emph{``convex polygonal set''}.  Note
that the word ``polygon'' is commonly used to refer to the convex hull
of a finite set of points in $\R^2$ which is a compact set (this is
not necessarily the case in algebraic geometry, e.g. Newton polygons).

\subsubsection{The polygon invariant}

Let $B_{\op{r}}:=\op{Int}(B)\setminus \{c_1,\,\ldots,\,c_{m_f}\}$,
which is precisely the set of regular values of $F$.  Given a sign
$\epsilon_i\in\{-1,+1\}$, let $\ell_i^{\epsilon_i}\subset\ell_i$ be
the vertical half line starting at $c_i$ and extending in the
direction of $\epsilon_i$~: upwards if $\epsilon_i=1$, downwards if
$\epsilon_i=-1$. Let $ \ell^{\vec\epsilon}:=
\bigcup_{i=1}^{m_f}\ell_i^{\epsilon_i}.  $ In Th.~3.8 in \cite{vungoc}
it was shown that:

\begin{theorem}
  For $\vec\epsilon\in\{-1,+1\}^{m_f}$ there exists a homeomorphism $f
  = f_\epsilon\colon B \to \R^2$, modulo a left composition by a
  transformation in $\mathfrak{I}$, such that $f|_{(B\setminus
    \ell^{\vec\epsilon})}$ is a diffeomorphism into the image of $f$,
  $\Delta:=f(B)$, which is a \emph{rational convex polygon},
  $f|_{(B_r\setminus \ell^{\vec\epsilon})}$ is affine (it sends the
  integral affine structure of $B_r$ to the standard structure of
  $\R^2$) and $f$ preserves $J$: i.e.  $
  f(x,\,y)=(x,\,f^{(2)}(x,\,y)).  $
\end{theorem}

The map $f$ satisfies further properties \cite{pelayovungoc1}, which
are relevant for the uniqueness proof.  In order to arrive at $\Delta$
one cuts $(J,\,H)(M) \subset \R^2$ along each of the vertical
half-lines $\ell_i^{\epsilon_i}$. Then the resulting image becomes
simply connected and thus there exists a global 2-torus action on the
preimage of this set. The polygon $\Delta$ is just the closure of the
image of a toric momentum map corresponding to this torus action.

We can see that this polygon is not unique. The choice of the ``cut
direction'' is encoded in the signs $\epsilon_j$, and there remains
some freedom for choosing the toric momentum map. Precisely, the
choices and the corresponding homeomorphisms $f$ are the following~:
\begin{itemize}
\item[(a)] {\em an initial set of action variables $f_0$ of the form
    $(J,\,K)$} near a regular Liouville torus in \cite[Step 2, pf. of
  Th.~3.8]{vungoc}.  If we choose $f_1$ instead of $f_0$, we get a
  polygon $\Delta'$ obtained by left composition with an element of
  $\mathfrak{I}$.  Similarly, if we choose $f_1$ instead of $f_0$, we
  obtain $f$ composed on the left with an element of $\mathfrak{I}$;
\item[(b)] {\em a tuple $\vec{\epsilon}$ of $1$ and $-1$}.  If we
  choose $\vec{\epsilon'}$ instead of $\vec{\epsilon}$ we get $
  \Delta'=t_{\vec{u}}(\Delta) $ with $u_i=(\epsilon_i-\epsilon'_i)/2$,
  by \cite[Proposition 4.1, expression (11)]{vungoc}.  Similarly instead of $f$
  we obtain $f'=t_{\vec{u}} \circ f$.
\end{itemize}

Once $f_0$ and $\vec\epsilon$ have been fixed as in (a) and (b),
respectively, then there exists a unique toric momentum map $\mu$ on
$M_r:=F^{-1}(\textup{Int}{B}\setminus(\bigcup \ell_j^{\epsilon_j}))$
which preserves the foliation $\mathcal{F}$, and coincides with
$f_0\circ F$ where they are both defined. Then, necessarily, the first
component of $\mu$ is $J$, and we have
% \begin{equation}
$\overline{\mu(M_r)}=\Delta$.
% \label{equ:mu-delta}
% \end{equation}
% \end{lemma}

We need now for our purposes to formalize choices (a) and (b) in a
single geometric object.  The details of how to do this have appeared
in \cite{pelayovungoc2}.  We will simply say that essentially this
object consists of the convex polygon itself together with a
collection of oriented cuts as in the figure below.  We call this
object a \emph{weighted polygon}, and denote it by
$\Delta_{\scriptop{w}}$.  The cuts are a collection of vertical lines
that go through the singularities.  The actual object is actually more
complex, as it is defined as an equivalence class of such polygon
considered as a part of a larger space of polygons on which several
groups act non\--trivially. These groups are $\{-1,\,1\}^{m_f}$ and
$\mathfrak{I}$.  The actual invariant is then denoted by
$[\Delta_{\scriptop{w}}]$.

\section{More symplectic invariants of semitoric
  systems} \label{sec:inv}

In \cite[Th.~6.2]{pelayovungoc1} the authors constructed, starting
from a given semitoric integrable system on a $4$\--manifold, a
collection of five symplectic invariants associated with it and proved
that these completely determine the integrable system up to global
isomorphisms of semitoric systems. Let $M_1$, $M_2$ be symplectic
$4$\--manifolds equipped with semitoric integrable systems
$(J_1,\,H_1)$ and $(J_2,\, H_2)$.  An \emph{isomorphism} between these
integrable systems is a symplectomorphism $ \varphi \colon M_1 \to M_2
$ such that $\varphi^*(J_2,\,H_2)=(J_1,\,f(J_1,\,H_1))$ for some
smooth function $f$ such that $\frac{\partial f}{\partial H_1}$
nowhere vanishes.

We recall the definition of the invariants that we assigned to a
semitoric integrable system in our previous paper
\cite{pelayovungoc1}, to which we refer to further details.  Then we
state the uniqueness theorem proved therein.

\subsection{Taylor series invariant}
\label{taylor:sec}

We assume that the critical fiber $ \mathcal{F}_m:=F^{-1}(c_i) $
contains only one critical point $m$, which according to Nguy{\^e}n
Ti{\^e}n Zung \cite{zung-I} is a generic condition, and let
$\mathcal{F}$ denote the associated singular foliation.  Moreover, we
will make for simplicity an even stronger generic assumption~: if $m$
is a focus-focus critical point for $F$, then $m$ is the unique
critical point of the level set $J^{-1}(J(m))$.  A semitoric system is
\emph{simple} if this generic assumption is satisfied.

These conditions imply that the values $J(m_1),\,\dots,\,J(m_{m_f})$
are pairwise distinct. We assume throughout the article that the
critical values $c_i$'s are \emph{ordered} by their $J$-values~:
$J(m_1)< J(m_2) < \cdots < J(m_{m_f})$.  Let $(S_i)^{\infty}$ be the a
formal power series expansion (in two variables with vanishing
constant term) corresponding to the integrable system given by $F$ at
the critical focus\--focus point $c_i$, see Theorem
\ref{theo:invariants}. We say that $(S_i)^\infty$ is the \emph{Taylor series invariant of
  $(M,\, \omega,\, (J,\,H))$ at the focus\--focus point $c_i$}.

\subsection{The Volume Invariant}
\label{volumesection}

Consider a focus\--focus critical point $m_i$ whose image by $(J,\,H)$
is $c_i$, and let $\Delta$ be a rational convex polygon corresponding
to the system $(M,\, \omega,\, (J,\,H))$.  If $\mu$ is a toric
momentum map for the system $(M, \, \omega,\,(J, \, H))$ corresponding
to $\Delta$, then the image $\mu(m_i)$ is a point in the interior of
$\Delta$, along the line $\ell_i$.  We proved in \cite{pelayovungoc1}
that the vertical distance
\begin{eqnarray} \label{height:eq} h_i:=\mu(m_i)-\min_{s \in \ell_i
    \cap \Delta} \pi_2(s)>0
\end{eqnarray}
is independent of the choice of momentum map $\mu$. Here $\pi_2 \colon
\R^2 \to \R$ is $\pi_2(x,\,y)=y$. The reasoning behind writing the
word ``volume'' in the name of this invariant is that it has the
following geometric interpretation: the singular manifold
$Y_i=J^{-1}(c_i)$ splits into $Y_i\cap\{H>H(m_i)\}$ and
$Y_i\cap\{H<H(m_i)\}$, and $h_i$ is the Liouville volume of
$Y_i\cap\{H<H(m_i)\}$.

\subsection{The Twisting\--Index Invariant}
\label{indexsection}
This is a subtle invariant of semitoric systems; it quantifies the
dynamical complexity of the system at a global level, while involving
the behavior near all of the focus\--focus singularities of the system
simultaneously.

The twisting-index expresses the fact that there is, in a neighborhood
of any focus\--focus critical point $c_i$, a \emph{privileged toric
  momentum map} $\nu$. This momentum map, in turn, is due to the
existence of a unique hyperbolic radial vector field in a neighborhood
of the focus-focus fiber. Therefore, one can view the twisting-index
as a dynamical invariant. Since any semitoric polygon defines a
(generalized) toric momentum map $\mu$, we will be able to define the
twisting-index as the integer $k_i\in \ZM$ such that
\[
\DD\mu = T^{k_i} \DD\nu.
\]
(Recall formula (\ref{equ:Tk}) for the formula of $T^{k_i}$.) We could
have equivalently defined the twisting-indices by comparing the
privileged momentum maps at different focus-focus points.

The precise definition of $k_i$ requires some care, which we explain
now.  Let $\Delta_{\scriptop{w}}$ be the weighted polygon associated
to $M$, which recall consists of the polygon $\Delta$ plus a
collection of oriented vertical lines $\ell_j$, where the orientation
of each line is given by $\pm 1$ signs $\epsilon_j$,
$j=1,\,\ldots,\,m_f$.

Let $\ell:=\ell_i^{\epsilon_i}\subset\R^2$ be the vertical
\emph{half\--line} starting at $c_i$ and pointing in the direction of
$\epsilon_i\, e_2$, where $e_1,\,e_2$ are the canonical basis vectors
of $\mathbb{R}^2$.

By Eliasson's theorem, there is a neighbourhood $W=W_i$ of the
focus\--focus critical point $m_i=F^{-1}(c_i)$, a local
symplectomorphism $\phi:(\R^4,\,0)\to W$, and a local diffeomorphism
$g$ of $(\R^2,\,0)$ such that $F \circ \phi= g\circ q$, where $q$ is
given by~\eqref{equ:cartan}.

Since $q_2\circ\phi^{-1}$ has a $2\pi$\--periodic Hamiltonian flow, it
is equal to $J$ in $W$, up to a sign.  Composing if necessary $\phi$
by $(x,\,\xi)\mapsto(-x,\,-\xi)$ one can assume that $q_1=J\circ\phi$ in
$W$, i.e. $g$ is of the form $g(q_1,\,q_2)=(q_1,\,g_2(q_1,\,q_2)).$
Upon composing $\phi$ with
$(x,\,y,\,\xi,\,\eta)\mapsto(-\xi,\,-\eta,\,x,\,y)$, which changes
$(q_1,\,q_2)$ into $(-q_1,\,q_2)$, one can assume that $
\frac{\partial g_2}{\partial q_2}(0)>0.  $ In particular, near the
origin $\ell$ is transformed by $g^{-1}$ into the positive imaginary
axis if $\epsilon_i=1$, or the negative imaginary axis if
$\epsilon_i=-1$.

\begin{figure}[htb]
  \begin{center}
    \includegraphics[width=0.4\linewidth]{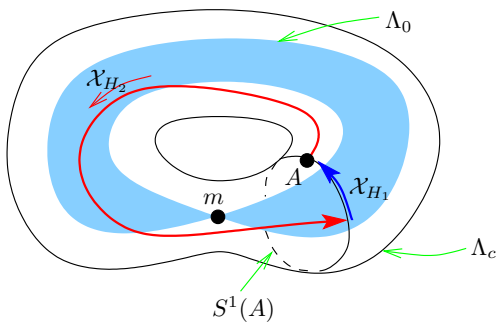}
    \caption{Singular foliation near the leaf $\Lambda_0=\mathcal{F}_m$, where
      $S^1(A)$ denotes the $S^1$\--orbit generated by $H_1=J$.}
  \end{center}
  \label{AF6}
\end{figure}

Let us now fix the origin of angular polar coordinates in $\R^2$ on
the \emph{positive} imaginary axis, let $V=F(W)$ and define
$\tilde{F}=(H_1,\,H_2)=g^{-1}\circ F$ on $F^{-1}(V)$ (notice that
$H_1=J$).

Recall that near any regular torus there exists a Hamiltonian vector
field $\mathcal{H}_p$, whose flow is $2\pi$\--periodic, defined by
  $$
  2\pi \mathcal{H}_p =
  (\tau_1\circ\tilde{F})\mathcal{H}_{H_1}+(\tau_2\circ\tilde{F})\mathcal{H}_J,
  $$
  where $\tau_1$ and $\tau_2$ are functions on $\R^2\setminus\{0\}$
  satisfying~\eqref{equ:sigma}, with $\sigma_2(0)>0$. In fact $\tau_1$
  is multivalued, but we determine it completely in polar coordinates
  with angle in $[0,\,2\pi)$ by requiring continuity in the angle
  variable and $\sigma_1(0)\in[0,\,2\pi)$. In case $\epsilon_i=1$,
  this defines $\mathcal{H}_p$ as a smooth vector field on
  $F^{-1}(V\setminus\ell)$.

  In case $\epsilon_i=-1$ we keep the same $\tau_1$\--value on the
  negative imaginary axis, but extend it by continuity in the angular
  interval $[\pi,\,3\pi)$. In this way $\mathcal{H}_p$ is again a
  smooth vector field on $F^{-1}(V\setminus\ell)$.
  
  Let $\mu$ be the generalized toric momentum map associated to
  $\Delta$. On $F^{-1}(V\setminus \ell$), $\mu$ is smooth, and its
  components $(\mu_1,\,\mu_2)=(J,\,\mu_2)$ are smooth Hamiltonians
  whose vector fields $(\mathcal{H}_J,\mathcal{H}_{\mu_2})$ are
  tangent to the foliation, have a $2\pi$-periodic flow, and are
  \emph{a.e.}  independent. Since the couple $(\mathcal{H}_J,
  \mathcal{H}_p)$ shares the same properties, there must be a matrix
  $A\in\textup{GL}(2,\Z)$ such that
  $(\mathcal{H}_J,\mathcal{H}_{\mu_2}) = A
  (\mathcal{H}_J,\mathcal{H}_p)$. This is equivalent to saying that
  there exists an integer $k_i\in\Z$ such that $ \mathcal{H}_{\mu_2} =
  k_i \mathcal{H}_{J} + \mathcal{H}_p.  $

  It was shown in \cite[Proposition~5.4]{pelayovungoc1} that $k_i$ is well
  defined, i.e. does not depend on choices.  The integer $k_i$ is
  called the \emph{twisting index of $\Delta_{\scriptop{w}}$ at the
    focus\--focus critical value $c_i$}.

  It was shown in \cite[Lemma~5.6]{pelayovungoc1} that there exists a
  unique smooth function $H_p$ on $F^{-1}(V\setminus \ell)$ with
  Hamiltonian vector field $\mathcal{H}_p$ and such that $\lim_{m\to
    m_i}H_p=0.$ The toric momentum map $\nu:=(J,\,H_p)$ is called {\em
    the privileged momentum map for $(J,\,H)$} around the
  focus\--focus value $c_i$.  If $k_i$ is the twisting index of $c_i$,
  one has $ \DD\mu =T^{k_i} \DD\nu $ on $F^{-1}(V)$. However, the
  twisting index does depend on the polygon $\Delta$. Thus, since we
  want to define an invariant of the initial semitoric system, we need
  to quotient out by the natural action of groups $G_{m_f}\times
  \mathfrak{I}$; because this is a rather technical task and we refer
  to \cite[p. 580]{pelayovungoc1} for details.

  It was shown in \cite[Proposition~5.8]{pelayovungoc1} that if two weighted
  polygons $\Delta_{\scriptop{w}}$ and $\Delta'_{\scriptop{weight}}$
  lie in the same $G_{m_f}$\--orbit, then the twisting indices
  $k_i,\,k'_i$ associated to $\Delta_{\scriptop{w}}$ and
  $\Delta'_{\scriptop{weight}}$ at their respective focus\--focus
  critical values $c_i,\,c'_i$ are equal.

  To a semitoric system we associate what we call the
  \emph{twisting\--index} invariant, which is nothing but the tuple
  $(\Delta_{\op{w}},\, {\bf k})$ consisting of the polygon $\Delta$
  labeled by the tuple twisting indices ${\bf
    k}=(k_j)_{j=1}^{m_f}$. Actually, as explained above, one needs to
  take into consideration the group actions of $G_{m_f}$ and
  $\mathfrak{I}$, so the twisting index invariant associated to the
  semitoric system is an equivalence class $ [(\Delta_{\op{w}}, \,
  {\bf k})] $ under a twisted action of $G_{m_f}$ and $\mathfrak{I}$.
  The formula for this action is long and we choose to not write it
  here, but details appear in \cite{pelayovungoc1}.

  \subsection{Example} In the case of the coupled spin oscillator the
  twisting index invariant does not appear because there is only one
  focus\--focus point. So in addition to the Taylor series invariant
  (of which as we said one can compute its linear approximation), the
  height invariant and the polygon invariant are easy to compute. They
  are explicitly given in Figure \ref{fig:spin} in the next section.

  \section{Global symplectic theory of semitoric
    systems} \label{sec:semitoric2}

  \subsection{First global result: uniqueness}

  The symplectic invariants constructed in \cite{pelayovungoc1}, for a
  given $4$\--dimensional semitoric integrable system, are the
  following: (i) \emph{the number of singularities invariant}: an
  integer $m_f$ counting the number of isolated singularities; (ii)
  \emph{the singularity type invariant}: a collection of $m_f$
  infinite formal Taylor series on two variables which classifies
  locally the type of each (focus\--focus) singularity; (iii)
  \emph{the polygon invariant}: the equivalence class of a weighted
  rational convex $\Delta_{\op{w}}$ consisting of a convex polygon
  $\Delta$ and the collection of vertical lines $\ell_j$ crossing it,
  where $\ell_j$ is oriented upwards or downwards depending on the
  sign of $\epsilon_j$, $j=1,\, \ldots,\, m_f$;
  \begin{figure}[htb]
    \begin{center}
      \includegraphics{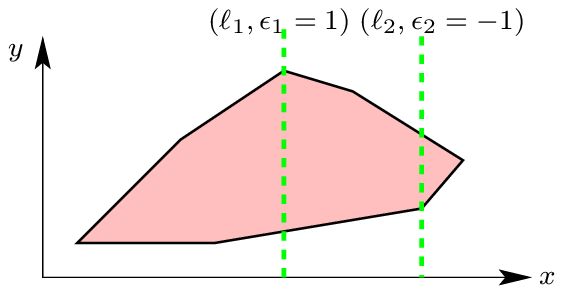}
      \caption{Weighted polygon $(\Delta, (\ell_1,\, \ell_2),\,
        (1,\,-1))$.}
    \end{center}
    \label{AF3}
  \end{figure}
  (iv) \emph{the volume invariant}: $m_f$ numbers measuring volumes of
  certain submanifolds at the singularities; (v) \emph{the twisting
    index invariant}: $m_f$ integers measuring how twisted the system
  is around singularities. This is a subtle invariant, which depends
  on the representative chosen in (iii).  Here, we write $m_f$ to
  emphasize that the singularities that $m_f$ counts are focus\--focus
  singularities.  We then proved:
  
 \begin{theorem}[Pelayo\--V\~ u Ng\d oc \cite{pelayovungoc1}]
   Two semitoric systems $(M_1,\, \omega_1,\,(J_1,\,H_1))$ and
   $(M_2,\, \omega_2,\, (J_2,\,H_2))$ are isomorphic if and only if
   they have the same invariants (i)--(v), where an \emph{isomorphism}
   is a symplectomorphism $ \varphi \colon M_1 \to M_2 $ such that
   $\varphi^*(J_2,\,H_2)=(J_1,\,f(J_1,\,H_1))$ for some smooth
   function $f$ such that $\frac{\partial f}{\partial H_1}$ nowhere
   vanishes.
 \end{theorem}

 \subsection{Second global result: existence}

 We have found that some restrictions on the symplectic invariants we
 have just defined must be imposed \cite{pelayovungoc2}.  Indeed, we
 call a ``semitoric list of ingredients'' the following collection of
 items:
 \begin{itemize}
 \item[(i)] An non\--negative integer $m_f$.
 \item[(ii)] An $m_{f}$\--tuple of formal Taylor series with vanishing
   constant term $((S_i)^{\infty})_{i=1}^{m_{f}} \in (\R[[X,\,Y]]_
   0)^{m_f}$.
 \item[(iii)]
   A Delzant semitoric polygon $[\Delta_{\scriptop{w}}]$ of complexity
   $m_f$ consisting of a polygon $\Delta$ and vertical lines $\ell_j$
   intersecting $\Delta$, each of which is oriented according to a
   sign $\epsilon_j=\pm1$;
 \item[(iv)] An $m_f$\--tuple of numbers ${\bf h}=(h_j)_{j=1}^{m_{f}}$
   such that $ 0 < h_j < \textup{length}(\Delta \cap \ell_i) $.
 \item[(v)] An equivalence class $[(\Delta_{\scriptop{w}},{\bf k})]$,
   where ${\bf k}=(k_j)_{j=1}^{m_f}$ is a collection of integers.
 \end{itemize}

 In the definition the term $\R[[X, \,Y]]$ refers to the algebra of
 real formal power series in two variables, and $\R[[X,\,Y]]_ 0$ is
 the subspace of such series with vanishing constant term, and first
 term $\sigma_1\,X+\sigma_2\,Y$ with $\sigma_2 \in [0,\,2\, \pi)$. For
 the definition of Delzant semitoric polygon, which is somewhat
 involved, see \cite[Sec. 4.2]{pelayovungoc2}. The main result of
 \cite{pelayovungoc2} is the following existence theorem:

  \begin{figure}[h]
    \centering
    \includegraphics[width=0.8\linewidth]{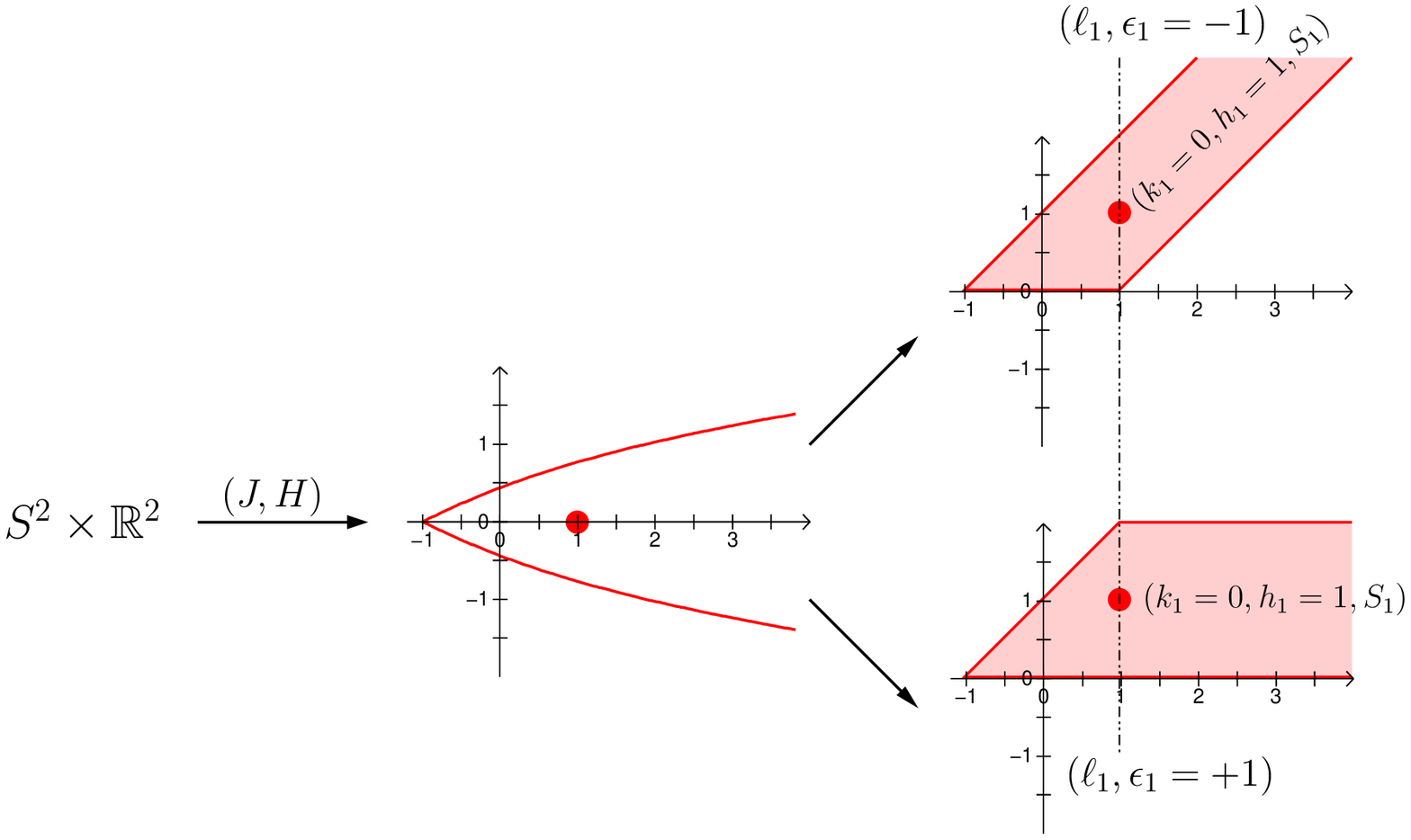}
    \caption{The coupled spin-oscillator example. The middle figure
      shows the image of the initial moment map $F=(J,\, H)$. Its
      boundary is the parameterized curve $(j(s)=\frac{s^2-3}{2s},
      h(s)=\pm\frac{s^2-1}{2s^{3/2}}), \,\,s \in[1,\infty).  $ The
      image is the connected component of the origin. The system is a
      simple semitoric system with one focus\--focus point whose image
      is $(1,\,0)$. The invariants are depicted on the right
      hand-side.  Since $m_f=1$, the class of generalized polygons for
      this system consists of two polygons.}
    \label{fig:spin}
  \end{figure}

  \begin{theorem}[Pelayo\--V\~ u Ng\d oc \cite{pelayovungoc2}]
    \label{class:thm}
    For each semitoric list of ingredients there exists a
    $4$\--dimensional simple semitoric integrable system with list of
    invariants equal to this list of ingredients.
  \end{theorem}

  The proof is involved, but the main idea of proof is simple.  We
  start with a representative of $[\Delta_{\scriptop{w}}]$ with all
  $\epsilon_{j}$'s equal to $+1$. The strategy is to construct the
  system locally and semiglobally around the singularities and around
  the regular parts, to then perform symplectic glueing in order to
  obtain a semitoric system by constructing a suitable singular torus
  fibration above $\Delta\subset\RM^2$.

  Rather subtle analytical issues appear when glueing, and one
  eventually ends up with a system given by a momentum map which is
  not smooth along the cuts $\ell_j$.  More concretely, first we
  construct a ``semitoric system'' over the part of the polygon away
  from the sets in the covering that contain the cuts $\ell_j^+$; then
  we attach to this ``semitoric system'' the focus\--focus fibrations
  i.e.  the models for the systems in a small neighborhood of the
  nodes (singularities).  We use a symplectic gluing theorem to do
  this glueing (cf. \cite{pelayovungoc2} for a statement/proof).

  Third, we continue to glue the local models in a small neighborhood
  of the cuts. The ``semitoric system'' is given by a proper toric map
  only in the preimage of the polygon away from the cuts. There is a
  analytically rather subtle issue near the cuts and one has to change
  the momentum map carefully to make it smooth while preserving the
  structure of the system up to isomorphisms.
 
  In the last step we prove that the system we have constructed has
  the right invariants. Here we have to appeal to the uniqueness
  theorem, as the equivalence class of the invariants may have shifted
  in the construction.

  \section{Some open problems} \label{sec:semitoric3}

  \subsection{Inverse spectral theory}

  Finding out how information from quantum completely integrable
  systems leads to information about classical systems is a
  fascinating ``inverse'' problem with very few precise results at
  this time.

  The \emph{symplectic} classification in terms of concrete invariants
  described in sections \ref{sec:semitoric}, \ref{sec:inv},
  \ref{sec:semitoric2} serves as a tool to quantize semitoric systems.
  In Delzant's theory, the image of the momentum map, for a toric
  completely integrable action, completely determines the system. In
  the quantum theory, the image of the momentum map is replaced by the
  joint spectrum. Can one determine the underlying classical system
  from the joint spectrum of the associated quantum system? In this
  vast, essentially unexplored program, one can ask the less ambitious
  but still spectacular question: given a quantum integrable system
  depending on a semiclassical parameter $\hbar$, and whose
  semiclassical limit $(J,H)$ is semitoric, does the knowledge of the
  semi-classical joint spectra for all values of $\hbar$ determine the
  underlying classical system $(J,H)$~?

 \begin{figure}[h]
   \centering
   \includegraphics[width=0.4\linewidth]{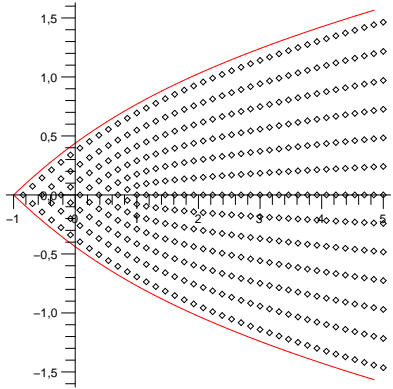}
   \caption{ Sections 4, 5 of the author's article \cite{example} are
     devoted to the spectral theory of the quantum coupled
     spin\--oscillator, and they are a first step towards proving this
     conjecture for spin\--oscillators. The unbounded operators
     $\hat{J}:=\op{Id} \otimes \Big( -\frac{\hbar^2}{2}
     \frac{\op{d}^2}{\op{d}u^2} +\frac{u^2}{2} \Big) + (\hat{z}
     \otimes \op{Id})$ and $\hat{H}=\frac{1}{2}(\hat{x}\otimes u +
     \hat{y} \otimes (\frac{\hbar}{\ii}\frac{\partial}{\partial u})$
     on the Hilbert space $ \mathcal{H} \otimes \op{L}^2(\R)\subset
     \op{L}^2(\R^2) \otimes \op{L}^2(\R)$ are self\--adjoint and
     commute, and they define the quantum spin\--oscillator. Their
     joint spectrum is depicted in the figure.}
   \label{fig:spectrumapprox2versionsmall}
 \end{figure}

 \begin{conj}
   A semitoric system $J,\, H$is determined up to symplectic
   equivalence by its semiclassical joint spectrum as $\hbar \to 0$.
   From any such spectrum one can construct explicitly the associated
   semitoric system, i.e.  the set of points in $\mathbb{R}^2$ where
   on the $x$\--basis we have the eigenvalues of $\hat{J}$, and on the
   vertical axis the eigenvalues of $\hat{H}$ restricted to the
   $\lambda$\--eigenspace of $\hat{J}$.
 \end{conj}

 The strategy to prove this is clear: given the joint spectrum, detect
 in it the symplectic invariants.  Once we have computed the
 symplectic invariants, we can symplectically recover the integrable
 system by \cite{pelayovungoc1, pelayovungoc2}, and hence the quantum
 system.  The authors have done this for the coupled\--spin oscillator
 \cite{example}. The method to recover the symplectic invariants from
 the joint spectrum combines microlocal analysis and Lie
 theory. Recovering the polygon invariant is probably the easiest and
 most pictorial procedure, as long as one stays on a heuristic
 level. Making the heuristic rigorous should be possible along the
 lines of the toric case explained in~\cite{san-book} and
 ~\cite{san-inverse}.

 \begin{figure}[h]
   \centering
   \includegraphics[width=0.3\textwidth]{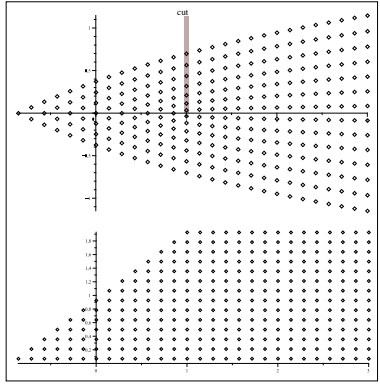}
   \caption{Recovering the polygon invariant. The top picture is the
     joint spectrum of $(\hat{J},\hat{H})$. In the bottom picture, we
     have developed the joint eigenvalues into a regular lattice. One
     can easily check on this illustration that the number of
     eigenvalues in each vertical line in the same in both pictures.}
   \label{fig:recover_polytope}
 \end{figure}
 The convex hull of the resulting set is a rational, convex polygonal
 set, depending on $\hbar$. Since the semiclassical affine structure
 is an $\hbar$-deformation of the classical affine structure, we see
 that, as $\hbar\to 0$, this polygonal set converges to the semitoric
 polygon invariant.

 \subsection{Mirror symmetry} \label{rmk:gross}

 When dealing with semitoric systems we are in a situation where the
 moment map $(J,H)$ is a ``torus fibration'' with singularities, and
 its base space becomes endowed with a singular integral affine
 structure. These same affine structures appear as a central
 ingredient in the work of Kontsevich and Soibelman \cite{KS}. These
 structures have been studied in the context of integrable systems (in
 particular by Nguy{\^e}n Ti{\^e}n Zung~\cite{zung-I}), but also
 became a central concept in the works by Symington \cite{s} and
 Symington\--Leung \cite{ls} in the context of symplectic geometry and
 topology, and by Gross\--Siebert, Casta{\~n}o\--Bernard,
 Castano\--Bernard\--Matessi \cite{grs1, grs2, grs3, grs4, castano-1,
   castano0, castano}, among others, in the context of mirror symmetry
 and algebraic geometry.  In fact the polygon invariant could have
 been expressed in terms of this affine structure.  It will be
 interesting to interpret the results of this paper in the context of
 mirror symmetry; at the least the classification of semitoric systems
 would give a large class of interesting examples. We hope to explore
 these ideas in the future.

 \subsection{Higher dimensions}
 
 It is natural to want to extend our $4$\--dimensional classification results to
 higher dimensions. We believe this is a very difficult
 problem in general. Physically, there is no reason for dimension 4 to
 be more or less relevant than higher dimensions.  The fact is that
 some of the results that our classification uses (primarily those of
 V\~ u Ng\d oc, but not exclusively) are $4$\--dimensional; but in
 principle there should be extensions to higher dimensions, as the
 proofs do not involve \emph{tools} that are specific to dimension 4.
   
 \subsection{Lagrange Top equations} The heavy top equations in body
 representation are known to be Hamiltonian on $\mathfrak{se}(3)^*$.
 These equations describe a classical Hamiltonian system with $2$
 degrees of freedom on the magnetic cotangent bundle $\op{T} S^2_{\|
   {\Gamma}_0\|} $. This two degrees of freedom system has a conserved
 integral but it does not have, generically, additional
 integrals. However, in the Lagrange case, one can find one additional
 integral which makes the system completely integrable.  It is
 classically known that the Lagrange heavy top is
 integrable. Moreover, one can check that it is semitoric, but it is
 given by a non\--proper momentum map.

 \begin{problem} \label{above} Develop the theory of semitoric systems
   $F=(J,\,H) \colon M \to \mathbb{R}^2$ when the component $J$
   generating a Hamiltonian circle action is not proper, but $F$ is proper.
 \end{problem}

 The authors' general theory does not cover the case stated in Problem
 \ref{above}, but many of techniques do extend at least to the case
 when $F$ is a proper map; we have been exploring this case in
 \cite{PeRaVu2010}.

\noindent
\\
{\'A}lvaro Pelayo \\
School of Mathematics\\
Institute for Advanced Study\\
Einstein Drive\\
Princeton, NJ 08540 USA.\\
 \\
\noindent
Washington University\\
Mathematics Department \\
One Brookings Drive, Campus Box 1146\\
St Louis, MO 63130-4899, USA.\\
{\em Website}: \url{http://www.math.wustl.edu/~apelayo/}\\
{\em E\--mail}: {apelayo@math.wustl.edu} and {apelayo@math.ias.edu}

\smallskip\noindent

\noindent
V\~u Ng\d oc San\\
Institut de Recherches Math{\'e}matiques de Rennes\\
Universit{\'e} de Rennes 1, Campus de Beaulieu\\
35042 Rennes cedex (France)\\
{\em Website}: \url{http://blogperso.univ-rennes1.fr/san.vu-ngoc/}\\
{\em E-mail:} {san.vu-ngoc@univ-rennes1.fr}

\end{document}